\pdfoutput=1
\RequirePackage{ifpdf}
\ifpdf 
\documentclass[pdftex]{sigma}
\else
\documentclass{sigma}
\fi

\usepackage[all]{xy}

\numberwithin{equation}{section}

\newtheorem{Theorem}{Theorem}[section]
\newtheorem{Lemma}[Theorem]{Lemma}
\newtheorem{Proposition}[Theorem]{Proposition}
\newtheorem{Corollary}[Theorem]{Corollary}
{\theoremstyle{definition}
\newtheorem{Definition}[Theorem]{Definition}
\newtheorem{Remark}[Theorem]{Remark}
}

\DeclareMathOperator{\ev}{ev}
\DeclareMathOperator{\ad}{ad}
\DeclareMathOperator{\ch}{ch}
\DeclareMathOperator{\supp}{supp}
\begin{document}

\allowdisplaybreaks

\renewcommand{\thefootnote}{$\star$}

\renewcommand{\PaperNumber}{047}

\FirstPageHeading

\ShortArticleName{Graded Limits of Minimal Af\/f\/inizations in Type $D$}

\ArticleName{Graded Limits of Minimal Af\/f\/inizations in Type $\boldsymbol{D}$\footnote{This paper is a~contribution to
the Special Issue on New Directions in Lie Theory. The full collection is available at \href{http://www.emis.de/journals/SIGMA/LieTheory2014.html}{http://www.emis.de/journals/SIGMA/LieTheory2014.html}}}

\Author{Katsuyuki NAOI}

\AuthorNameForHeading{K.~Naoi}

\Address{Institute of Engineering, Tokyo University of Agriculture and Technology,\\
3-8-1 Harumi-cho, Fuchu-shi, Tokyo, Japan} 
\Email{\href{naoik@cc.tuat.ac.jp}{naoik@cc.tuat.ac.jp}}

\ArticleDates{Received October 30, 2013, in f\/inal form April 14, 2014; Published online April 20, 2014}

\vspace{-1.5mm}

\Abstract{We study the graded limits of minimal af\/f\/inizations over a~quantum loop algebra of type $D$ in the regular
case.
We show that the graded limits are isomorphic to multiple generalizations of Demazure modules, and also give their
def\/ining relations.
As a~corollary we obtain a~character formula for the minimal af\/f\/inizations in terms of Demazure operators, and
a~multiplicity formula for a~special class of the minimal af\/f\/inizations.}

\Keywords{minimal af\/f\/inizations; quantum af\/f\/ine algebras; current algebras}

\Classification{17B37; 17B10}

\renewcommand{\thefootnote}{\arabic{footnote}}
\setcounter{footnote}{0}

\vspace{-3.5mm}

\section{Introduction}

\looseness=-1
Let ${\mathfrak g}$ be a~complex simple Lie algebra, $\mathbf{L}{\mathfrak g} = {\mathfrak g} \otimes
\mathbb{C}[t,t^{-1}]$ the associated loop algebra, and $U_q(\mathbf{L}{\mathfrak g})$ the quantum loop algebra.
In~\cite{MR1367675}, Chari introduced an important class of f\/inite-dimensional simple $U_q(\mathbf{L}{\mathfrak
g})$-modules called minimal af\/f\/inizations.
For a~simple $U_q({\mathfrak g})$-module $V$, we say a~simple $U_q(\mathbf{L}{\mathfrak g})$-module $\widehat{V}$ is an
af\/f\/inization of $V$ if the highest weight of $\widehat{V}$ is equal to that of~$V$.
One can def\/ine a~partial ordering on the equivalence classes (the isomorphism classes as a~$U_q({\mathfrak g})$-module)
of af\/f\/inizations of~$V$, and modules belonging to minimal classes are called minimal af\/f\/i\-ni\-za\-tions (a precise def\/inition
is given in Section~\ref{subsection:minimal_aff}).
For example, a~Kirillov--Reshetikhin module is a~minimal af\/f\/inization whose highest weight is a~multiple of
a~fundamental weight.
Minimal af\/f\/inizations have been the subjects of many articles in the recent years.
See~\cite{MR1347873,MR1402568,MR2342293,MR2587436,MY,N} for instance.
For the original motivations of considering minimal af\/f\/inizations, see~\cite[Introduction]{MR1367675}.
Given a~minimal af\/f\/inization, one can consider its classical limit.
By restricting it to the current algebra ${\mathfrak g}[t] = {\mathfrak g} \otimes \mathbb{C}[t]$ and taking
a~pull-back, a~graded ${\mathfrak g}[t]$-module called \textit{graded limit} is obtained.
Graded limits are quite important for the study of minimal af\/f\/inizations since the $U_q({\mathfrak g})$-module structure
of a~minimal af\/f\/inization is completely determined by the $U({\mathfrak g})$-module structure of its graded limit.

Graded limits of minimal af\/f\/inizations were f\/irst studied in~\cite{MR1836791,MR2238884} in the case of
Kirillov--Reshetikhin modules, and subsequently the general ones were studied in~\cite{MR2587436}.
In that paper, Moura presented several conjectures for the graded limits of minimal af\/f\/inizations in general types, and
partially proved them.
Graded limits of minimal af\/f\/inizations in type $ABC$ were further studied in~\cite{N}.
In that paper the author proved that the graded limit of a~minimal af\/f\/inization in these types is isomorphic to
a~certain ${\mathfrak g}[t]$-module $D(w_\circ\xi_1,\ldots,w_\circ\xi_n)$.
Here $w_\circ$ is the longest element of the Weyl group of ${\mathfrak g}$, $\xi_j$ are certain weights of the af\/f\/ine
Lie algebra $\widehat{{\mathfrak g}}$ which are ${\mathfrak g}$-dominant, and $D(w_\circ \xi_1,\ldots,w_\circ \xi_n)$ is
a~${\mathfrak g}[t]$-submodule of a~tensor product of simple highest weight $\widehat{{\mathfrak g}}$-modules, which is
generated by the tensor product $v_{w_\circ \xi_1} \otimes \dots \otimes v_{w_\circ \xi_n}$ of the extremal weight
vectors with weights $w_\circ\xi_j$.
As a~corollary of this, a~character formula for minimal af\/f\/inizations was given in terms of Demazure operators.
In addition, the def\/ining relations of graded limits conjectured in~\cite{MR2587436} were also proved.

In type $ABC$ a~minimal af\/f\/inization with a~f\/ixed highest weight is unique up to equivalence, and the graded limit of
a~minimal af\/f\/inization depends only on the equivalence class.
As a~consequence, the module $D(w_\circ \xi_1,\ldots, w_\circ \xi_n)$ can be determined from the highest weight only.
(If the highest weight is $\lambda= \sum\limits_{1 \le i \le n} \lambda_i\varpi_i$ where $\varpi_i$ are fundamental
weights, then $\xi_j$ are roughly equal to $\lambda_i(\varpi_i + a_i\Lambda_0)$ where $\Lambda_0$ is the fundamental
weight of $\widehat{{\mathfrak g}}$ associated with the distinguished node $0$, and $a_i = 1$ if the simple root $\alpha_i$
is long and $a_i = 1/2$ otherwise.
For the more precise statement, see~\cite{N}.)

In contrast to this, in type $D$ there are nonequivalent minimal af\/f\/inizations with the same highest weights.
It was proved in~\cite{MR1402568}, however, that even in type $D$ if the given highest weight satisf\/ies some mild
condition (see Section~\ref{subsection:minimal_aff}), then there are at most $3$ equivalence classes of minimal
af\/f\/inizations with the given highest weight.
We say a~minimal af\/f\/inization is \textit{regular} if its highest weight satisf\/ies this condition.
The purpose of this paper is to study the graded limits of regular minimal af\/f\/inizations of type~$D$ using the methods
in~\cite{N}.

In the sequel we assume that ${\mathfrak g}$ is of type $D_n$.
Let $\boldsymbol{\pi}$ be Drinfeld polynomials and assume that the simple $U_q(\mathbf{L}{\mathfrak g})$-module $L_q(\boldsymbol{\pi})$
associated with $\boldsymbol{\pi}$ is a~regular minimal af\/f\/inization.
Then in a~certain way we can associate with $L_q(\boldsymbol{\pi})$ a~vertex $s\in \{1,n-1,n\}$ of the Dynkin diagram of
${\mathfrak g}$ (see Section~\ref{subsection:minimal_aff}).
In the case where the number of equivalence classes are exactly~$3$, this $s$ parameterizes the equivalence class of
$L_q(\boldsymbol{\pi})$.
In this paper we show that there exists a~sequence $\xi_1^{(s)},\ldots,\xi_n^{(s)}$ of ${\mathfrak g}$-dominant
$\widehat{{\mathfrak g}}$-weights such that the graded limit $L(\boldsymbol{\pi})$ of $L_q(\boldsymbol{\pi})$ is isomorphic to
$D\big(w_\circ \xi_1^{(s)},\ldots, w_\circ \xi_n^{(s)}\big)$ (Theorem~\ref{Thm:Main2}).
Here~$\xi_j^{(s)}$ depends not only on the highest weight of~$L_q(\boldsymbol{\pi})$ but also~$s$, and the correspondence is
less straightforward compared with the case of type~$ABC$ (see Section~\ref{subsection:main_theorems} for the precise
statement).
As a~consequence, we give a~character formula for~$L_q(\boldsymbol{\pi})$ in terms of Demazure operators
(Corollary~\ref{Cor:characters}).
We also prove the def\/ining relations of the graded limits $L(\boldsymbol{\pi})$ conjectured in~\cite{MR2587436}
(Theorem~\ref{Thm:Theorem1}), which also depends not only on the highest weight but also $s$.

Recently Sam proved in~\cite{Sam} some combinatorial identity in type $BCD$, and gave a~multiplicity formula for minimal
af\/f\/inizations in type~$BC$ using the identity and results in~\cite{MR2763623} and~\cite{N}.
By applying the identity of type~$D$ to our results, we also obtain a~similar multiplicity formula for a~special class
of minimal af\/f\/inizations in type~$D$, which gives multiplicities in terms of the simple Lie algebra of type~$C$
(Corollary~\ref{Cor:Sam}).

The proofs of most results are similar to those in~\cite{N} and are in some respects even simpler since the type $D$ is
simply laced.
For example we do not need the theory of $q$-characters, which was essentially needed in \textit{loc.~cit.}

The organization of the paper is as follows.
In Section~\ref{Section2}, we give preliminary def\/initions and basic results.
In particular, we recall the def\/inition of the modules $D(\xi_1,\ldots,\xi_p)$, the classif\/ication of regular minimal
af\/f\/inizations of type~$D$, and the def\/inition of graded limits.
In Section~\ref{Section3} we state Theorems~\ref{Thm:Main2} and~\ref{Thm:Theorem1}, and discuss some of their corollaries.
The proofs of Theorems~\ref{Thm:Main2} and~\ref{Thm:Theorem1} is given in Section~\ref{section:Proof}.

\vspace{-1mm}

\section{Preliminaries}\label{Section2}
\vspace{-1mm}

\subsection[Simple Lie algebra of type $D_n$]{Simple Lie algebra of type $\boldsymbol{D_n}$}

Let $\widehat{I} = \{0,1,\ldots,n\}$ and $\widehat{C} = (c_{ij})_{i,j \in \widehat{I}}$ be the Cartan matrix of type $D_n^{(1)}$
whose Dynkin diagram is as follows:\vspace{-2mm}
\begin{gather}
\label{eq:Dynkin}
\begin{split}
\xygraph{ \circ ([]!{+(0,-.3)} {1}) - [r] \circ ([]!{+(0,-.3)} {2}) (- [u] \circ []!{+(0,.3)} {0}) - [r] \circ
([]!{+(0,-.3)} {3}) - [r] \dots - [r] \circ ([]!{+(0,-.3)} {n-3}) - [r] \circ ([]!{+(0,-.3)} {n - 2}) (- [u] \circ
[]!{+(0,.3)} {n})- [r] \circ ([]!{+(0,-.3)} {n-1}) }
\end{split}
\end{gather}
Let $J \subseteq \widehat{I}$ be a~subset.
In this paper, by abuse of notation, we sometimes denote by $J$ the subdiagram of~\eqref{eq:Dynkin} whose vertices are~$J$.

Let $I = \widehat{I} \setminus \{0\}$, $C = (c_{ij})_{i,j \in I}$ be the Cartan matrix of type~$D_n$, and ${\mathfrak g}$
the complex simple Lie algebra associated with~$C$.
Let~${\mathfrak h}$ be a~Cartan subalgebra and~${\mathfrak b}$ a~Borel subalgebra containing~${\mathfrak h}$.
Denote by $\Delta$ the root system and by~$\Delta_+$ the set of positive roots, and let $\theta \in \Delta_+$ be the
highest root.
Let $\alpha_i$ and $\varpi_i$ ($i\in I$) be the simple roots and fundamental weights respectively, and set $\varpi_0 =
0$ for convenience.
Let~$P$ be the weight lattice and $P^+$ the set of dominant integral weights.
Let~$W$ denote the Weyl group with simple ref\/lections~$s_i$ ($i \in I$), and~$w_\circ \in W$ the longest element.

For each $\alpha \in \Delta$ denote by ${\mathfrak g}_\alpha$ the corresponding root space, and f\/ix nonzero elements
$e_\alpha \in {\mathfrak g}_\alpha$, $f_\alpha \in {\mathfrak g}_{-\alpha}$ and $\alpha^\vee \in {\mathfrak h}$ such
that
\begin{gather*}
[e_\alpha, f_\alpha] = \alpha^\vee,
\qquad
[\alpha^\vee, e_\alpha] = 2 e_\alpha,
\qquad
[\alpha^\vee, f_\alpha] = -2 f_\alpha.
\end{gather*}
We also use the notation $e_i = e_{\alpha_i}$, $f_i = f_{\alpha_i}$ for $i \in I$.
Set ${\mathfrak n}_{\pm} = \bigoplus_{\alpha \in \Delta_+} {\mathfrak g}_{\pm \alpha}$.
For a~subset $J \subseteq I$, denote by ${\mathfrak g}_J \subseteq {\mathfrak g}$ the semisimple Lie subalgebra
corresponding to $J$, and let ${\mathfrak h}_J = \sum\limits_{i \in J} \mathbb{C} \alpha_i^{\vee} \subseteq {\mathfrak h}$.

\subsection[Af\/f\/ine Lie algebra of type $D_n^{(1)}$]{Af\/f\/ine Lie algebra of type $\boldsymbol{D_n^{(1)}}$} 

Let $\widehat{{\mathfrak g}} = {\mathfrak g} \otimes \mathbb{C}[t,t^{-1}] \oplus \mathbb{C} K \oplus \mathbb{C} d$ be the
af\/f\/ine Lie algebra with Cartan matrix $\widehat{C}$, where $K$ is the canonical central element and $d$ is the degree
operator.
Naturally ${\mathfrak g}$ is regarded as a~Lie subalgebra of $\widehat{{\mathfrak g}}$.
Def\/ine a~Cartan subalgebra $\widehat{{\mathfrak h}}$ and a~Borel subalgebra $\widehat{{\mathfrak b}}$ as follows:
\begin{gather*}
\widehat{{\mathfrak h}} = {\mathfrak h} \oplus \mathbb{C} K \oplus \mathbb{C} d,
\qquad
\widehat{{\mathfrak b}} = \widehat{{\mathfrak h}} \oplus {\mathfrak n}_+ \oplus {\mathfrak g} \otimes t\mathbb{C}[t].
\end{gather*}
Set $\widehat{{\mathfrak n}}_+ = {\mathfrak n}_+ \oplus {\mathfrak g} \otimes t\mathbb{C}[t]$.
We often consider ${\mathfrak h}^*$ as a~subspace of $\widehat{{\mathfrak h}}^*$ by setting $\langle K, \lambda \rangle =
\langle d, \lambda\rangle = 0$ for $\lambda \in {\mathfrak h}^*$.
Let $\widehat{\Delta}$ be the root system of $\widehat{{\mathfrak g}}$, $\widehat{\Delta}_+$ the set of positive roots,
$\widehat{\Delta}^{\mathrm{re}}$ the set of real roots and $\widehat{\Delta}^{\mathrm{re}}_+ = \widehat{\Delta}_+ \cap
\widehat{\Delta}^{\mathrm{re}}$.
Set $\alpha_0 = \delta - \theta$, $e_0 = f_\theta \otimes t$, $f_0 = e_\theta \otimes t^{-1}$ and $\alpha_0^\vee = K -
\theta^\vee$.

Denote by $\Lambda_0 \in \widehat{{\mathfrak h}}^*$ the unique element satisfying $\langle K, \Lambda_0 \rangle = 1$ and
$\langle {\mathfrak h}, \Lambda_0 \rangle = \langle d, \Lambda_0 \rangle = 0$, and def\/ine $\widehat{P}, \widehat{P}^+ \subseteq
\widehat{{\mathfrak h}}^*$ by
\begin{gather*}
\widehat{P} = P \oplus \mathbb{Z} \Lambda_0 \oplus \mathbb{C} \delta
\qquad
\text{and}
\qquad
\widehat{P}^+ = \big\{ \xi \in \widehat{P} \,|\, \langle \alpha_i^\vee, \xi \rangle \ge 0
\text{ for all }
i \in \widehat{I}\,\big\}.
\end{gather*}
Let $\widehat{W}$ denote the Weyl group of $\widehat{{\mathfrak g}}$ with simple ref\/lections $s_i$ ($i \in
\widehat{I}$).
We regard $W$ naturally as a~subgroup of $\widehat{W}$.
Let $\ell\colon \widehat{W} \to \mathbb{Z}_{\ge 0}$ be the length function.
Let $(\,,\,)$ be the unique non-degenerate $\widehat{W}$-invariant symmetric bilinear form on $\widehat{{\mathfrak h}}^*$
satisfying
\begin{gather*}
(\alpha, \alpha) = 2
\quad
\text{for}
\qquad
\alpha \in \widehat{\Delta}^{\mathrm{re}},
\qquad
({\mathfrak h}^*, \delta) = ({\mathfrak h}^*, \Lambda_0) = (\Lambda_0,\Lambda_0)= 0
\qquad
\text{and}
\qquad
(\delta, \Lambda_0) = 1.
\end{gather*}

Let $\Sigma$ be the group of Dynkin diagram automorphisms of $\widehat{{\mathfrak g}}$, which naturally acts on~$\widehat{{\mathfrak h}}^*$ and~$\widehat{{\mathfrak g}}$, and $\widetilde{W}$ the subgroup of ${\rm GL}(\widehat{{\mathfrak h}}^*)$
generated by $\widehat{W}$ and $\Sigma$.
Note that we have $\widetilde{W} = \Sigma \ltimes \widehat{W}$.
The length function $\ell$ is extended on $\widetilde{W}$ by setting $\ell(\tau w) = \ell(w)$ for $\tau \in \Sigma$, $w
\in \widehat{W}$.

Denote by $V(\lambda)$ for $\lambda \in P^+$ the simple ${\mathfrak g}$-module with highest weight $\lambda$, and by
$\widehat{V}(\Lambda)$ for $\Lambda\in \widehat{P}^+$ the simple highest weight $\widehat{{\mathfrak g}}$-module with highest
weight $\Lambda$.
For a~f\/inite-dimensional semisimple ${\mathfrak h}$-module (resp.~$\widehat{{\mathfrak h}}$-module) $M$ we denote by
$\ch_{{\mathfrak h}} M \in \mathbb{Z}[{\mathfrak h}^*]$ (resp.~$\ch_{\widehat{{\mathfrak h}}} M \in
\mathbb{Z}[\widehat{{\mathfrak h}}^*]$) its character with respect to ${\mathfrak h}$ (resp.~$\widehat{{\mathfrak h}}$).
We will omit the subscript ${\mathfrak h}$ or $\widehat{{\mathfrak h}}$ when it is obvious from the context.

\subsection{Loop algebras and current algebras}

Given a~Lie algebra $\mathfrak{a}$, its \textit{loop algebra} $\mathbf{L}{\mathfrak a}$ is def\/ined as the tensor product
${\mathfrak a} \otimes \mathbb{C}[t,t^{-1}]$ with the Lie algebra structure given by $[x \otimes f, y \otimes g] =
[x,y]\otimes fg$.
Let ${\mathfrak a}[t]$ and $t^k {\mathfrak a}[t]$ for $k \in \mathbb{Z}_{> 0}$ denote the Lie subalgebras ${\mathfrak a}
\otimes \mathbb{C}[t]$ and ${\mathfrak a} \otimes t^k \mathbb{C}[t]$ respectively.
The Lie algebra ${\mathfrak a}[t]$ is called the \textit{current algebra} associated with ${\mathfrak a}$.

For $a \in \mathbb{C}^\times$, let $\ev_a\colon \mathbf{L}{\mathfrak g} \to {\mathfrak g}$ denote the
\textit{evaluation map} def\/ined by $\ev_a(x \otimes f) = f(a)x$, and let $V(\lambda,a)$ for $\lambda \in P^+$ be
the \textit{evaluation module} which is the simple $\mathbf{L}{\mathfrak g}$-module def\/ined by the pull-back of
$V(\lambda)$ with respect to $\ev_a$.
An evaluation module for ${\mathfrak g}[t]$ is def\/ined similarly and is denoted by $V(\lambda,a)$ ($\lambda \in P^+$,
$a~\in \mathbb{C}$).

\subsection[$\widehat{\mathfrak b}$-submodules $D(\xi_1,\ldots,\xi_p)$]{$\boldsymbol{\widehat{\mathfrak b}}$-submodules $\boldsymbol{D(\xi_1,\ldots,\xi_p)}$}
\label{subsection:b-module}

Let $\xi_1,\ldots,\xi_p$ be a~sequence of elements belonging to the Weyl group orbits $\widehat{W} (\widehat{P}^+)$ of dominant
integral weights of $\widehat{{\mathfrak g}}$.
We def\/ine a~$\widehat{{\mathfrak b}}$-module $D(\xi_1,\ldots,\xi_p)$ as follows.
For each $1 \le j \le p$ let $\Lambda^j \in \widehat{P}^+$ be the unique element satisfying $\xi_j \in \widehat{W} \Lambda^j$,
and take a~nonzero vector $v_{\xi_j}$ in the $1$-dimensional weight space $\widehat{V}(\Lambda^j)_{\xi_j}$.
Then def\/ine
\begin{gather*}
D(\xi_1, \ldots,\xi_p) = U(\widehat{{\mathfrak b}})(v_{\xi_1} \otimes \dots \otimes v_{\xi_{p}}) \subseteq
\widehat{V}(\Lambda^1) \otimes \dots \otimes \widehat{V}(\Lambda^{p}).
\end{gather*}
If $(\alpha_i,\xi_j) \le 0$ for all $i \in I$ and $1 \le j \le p$, then $D(\xi_1, \ldots,\xi_p)$ can be regarded as
a~${\mathfrak g}[t] \oplus \mathbb{C} K \oplus \mathbb{C} d$-module and in particular a~${\mathfrak g}[t]$-module.

Some of $D(\xi_1,\ldots,\xi_p)$ are realized in a~dif\/ferent way.
To introduce this, we prepare some notation.
Assume that $V$ is a~$\widehat{{\mathfrak g}}$-module and $D$ is a~$\widehat{{\mathfrak b}}$-submodule of $V$.
For $\tau \in \Sigma$, we denote by~$F_\tau V$ the pull-back $(\tau^{-1})^*V$ with respect to the Lie algebra
automorphism $\tau^{-1}$ on $\widehat{{\mathfrak g}}$, and def\/ine a~$\widehat{{\mathfrak b}}$-submodule $F_\tau D \subseteq
F_\tau V$ in the obvious way.
It is easily proved that
\begin{gather*}
F_\tau D(\xi_1,\ldots,\xi_p) \cong D(\tau\xi_1,\ldots,\tau\xi_p).
\end{gather*}
For $i \in \widehat{I}$ let $\widehat{{\mathfrak p}}_i$ denote the parabolic subalgebra $\widehat{{\mathfrak b}} \oplus \mathbb{C}
f_i \subseteq \widehat{{\mathfrak g}}$, and set $F_iD= U(\widehat{{\mathfrak p}}_i)D \subseteq V$ to be the $\widehat{{\mathfrak
p}}_i$-submodule generated by $D$.
Finally we def\/ine $F_w D$ for $w \in \widetilde{W}$ as follows: let $\tau \in \Sigma$ and $w' \in \widehat{W}$ be the
elements such that $w=\tau w'$, and choose a~reduced expression $w'=s_{i_1} \cdots s_{i_k}$.
Then we set
\begin{gather*}
F_{w}D = F_\tau F_{i_1}\cdots F_{i_k} D \subseteq F_\tau V.
\end{gather*}

\begin{Proposition}[\protect{\cite[Proposition 2.7]{N}}]
\label{Prop:construction}
Let $\Lambda^1, \ldots, \Lambda^p$ be a~sequence of elements of $\widehat{P}^+$, and $w_{1}, \ldots, w_{p}$ a~sequence of
elements of $\widetilde{W}$.
We write $w_{[r,s]}= w_{r} w_{r+1}\cdots w_s$ for $r \le s$, and assume that $\ell(w_{[1,p]}) = \sum\limits_{j=1}^p
\ell(w_{j})$.
Then we have
\begin{gather}
D\big(w_{[1,1]} \Lambda^1, w_{[1,2]} \Lambda^2, \ldots,w_{[1,p-1]}\Lambda^{p-1}, w_{[1, p]}\Lambda^p\big)
\nonumber
\\
\qquad{}
\cong F_{w_{1}}\big(D\big(\Lambda^1\big) \otimes F_{w_{2}}\big(D\big(\Lambda^2\big)\otimes \cdots \otimes
F_{w_{p-1}}\big(D\big(\Lambda^{p-1}\big) \otimes F_{w_{p}} D\big(\Lambda^p\big)\big)\cdots \big)\big).
\label{eq:Demazure_realizatin}
\end{gather}
\end{Proposition}

Let $\mathcal{D}_i$ for $i \in \widehat{I}$ be a~linear operator on $\mathbb{Z}[\widehat{P}]$ def\/ined by
\begin{gather*}
\mathcal{D}_i(f) = \frac{f - e^{-\alpha_i}s_i(f)}{1 - e^{-\alpha_i}},
\end{gather*}
where $e^\lambda$ ($\lambda \in \widehat{P}$) are the generators of $\mathbb{Z}[\widehat{P}]$.
For $w \in \widehat{W}$ with a~reduced expression $w= s_{i_1}\cdots s_{i_k}$, we set $\mathcal{D}_w =
\mathcal{D}_{i_1}\cdots\mathcal{D}_{i_k}$.
If $w \in \widetilde{W}$ and $w=\tau w'$ ($\tau \in \Sigma$, $w' \in \widehat{W}$), we set $\mathcal{D}_w =
\tau\mathcal{D}_{w'}$.
The operator $\mathcal{D}_w$ is called a~\textit{Demazure operator}.
The character of the right-hand side of~\eqref{eq:Demazure_realizatin} is expressed using Demazure operators
by~\cite[Theorem~5]{MR1887117}, and as a~consequence we have the following (see also~\cite[Corollary~2.8]{N}).

\begin{Proposition}
\label{Prop:character}
Let $\Lambda^j \in \widehat{P}^+$ and $w_j \in \widetilde{W}$ $(1 \le j \le p)$ be as in
Proposition~{\rm \ref{Prop:construction}}.
Then we have
\begin{gather*}
\ch_{\widehat{{\mathfrak h}}} D\big(w_{[1,1]} \Lambda^1,
w_{[1,2]} \Lambda^2, \ldots, w_{[1,p-1]}\Lambda^{p-1},
w_{[1, p]}\Lambda^p\big)
\\
\qquad{}
= \mathcal{D}_{w_1}\big(e^{\Lambda^1}\cdot \mathcal{D}_{w_2}\big(e^{\Lambda^2}\cdots
\mathcal{D}_{w_{p-1}}\big(e^{\Lambda^{p-1}}\cdot \mathcal{D}_{w_p} \big(e^{\Lambda^p}\big)\big)\cdots\big)\big).
\end{gather*}
\end{Proposition}

\subsection{Quantum loop algebras and their representations}

The quantum loop algebra $U_q(\mathbf{L}{\mathfrak g})$ is a~$\mathbb{C}(q)$-algebra generated by $x_{i,r}^{\pm}$,
$k_i^{\pm 1}$ and $h_{i,m}$ ($i \in I$, $r \in \mathbb{Z}$, $m \in \mathbb{Z} \setminus \{0\}$) subject to certain
relations (see, e.g.,~\cite[Section 12.2]{MR1300632}).
$U_q(\mathbf{L}{\mathfrak g})$ has a~Hopf algebra structure~\cite{MR1300632,MR1227098}.
In particular if $V$ and $W$ are $U_q(\mathbf{L}{\mathfrak g})$-modules then $V\otimes W$ and $V^*$ are also
$U_q(\mathbf{L}{\mathfrak g})$-modules, and we have $(V \otimes W)^* \cong W^* \otimes V^*$.

Denote by $U_q(\mathbf{L}{\mathfrak n}_{\pm})$ and $U_q(\mathbf{L}{\mathfrak h})$ the subalgebras of
$U_q(\mathbf{L}{\mathfrak g})$ generated by $\{x_{i,r}^\pm \,|\, i \in I, r \in \mathbb{Z} \}$ and $\big\{k^{\pm
1}_i,h_{i,m} \,|\, i \in I, m \in \mathbb{Z} \setminus \{0\} \big\}$ respectively.
Denote by $U_q({\mathfrak g})$ the subalgebra generated by $\{x_{i,0}^{\pm}, k^{\pm 1}_i \,|\, i \in I\}$, which is
isomorphic to the quantized enveloping algebra associated with ${\mathfrak g}$.
For a~subset $J \subseteq I$, let $U_q(\mathbf{L}{\mathfrak g}_J)$ denote the subalgebra generated by
$\big\{x_{i,r}^{\pm},k_i^{\pm 1}, h_{i,m}\,|\, i\in J$, $r \in \mathbb{Z}, \, m \in \mathbb{Z} \setminus \{ 0\} \big\}$.

We recall basic results on f\/inite-dimensional $U_q({\mathfrak g})$- and $U_q(\mathbf{L}{\mathfrak g})$-modules.
Note that in the present paper we assume that ${\mathfrak g}$ is of type $D$, and when ${\mathfrak g}$ is non-simply
laced some of indeterminates $q$ appearing below should be replaced by $q_i = q^{d_i}$ with suitable $d_i \in
\mathbb{Z}_{>0}$.

A $U_q({\mathfrak g})$-module (or $U_q(\mathbf{L}{\mathfrak g})$-module) $V$ is said to be \textit{of type $1$} if $V$  
satisf\/ies
\begin{gather*}
V = \bigoplus_{\lambda \in P} V_\lambda,
\qquad
V_\lambda = \big\{ v \in V \,|\, k_iv = q^{\langle \alpha_i^\vee, \lambda \rangle}v\big\}.
\end{gather*}
In this article we will only consider modules of type~$1$.
For a~f\/inite-dimensional module~$V$ of type~$1$, we set $\ch V = \sum\limits_{\lambda \in P} e^\lambda\dim
V_\lambda \in \mathbb{Z}[P]$.
For $\lambda \in P^+$ we denote by $V_q(\lambda)$ the f\/inite-dimensional simple $U_q({\mathfrak g})$-module of type $1$
with highest weight $\lambda$.
The category of f\/inite-dimensional $U_q({\mathfrak g})$-modules of type $1$ is semisimple, and every simple object is
isomorphic to~$V_q(\lambda)$ for some $\lambda \in P^+$.

We say that a~$U_q(\mathbf{L}{\mathfrak g})$-module $V$ is \textit{highest $\ell$-weight} with highest $\ell$-weight  
vector $v$ and highest $\ell$-weight $\big(\gamma^+_i(u), \gamma^-_i(u)\big)_{i\in I} \in \big(\mathbb{C}(q)[[u]] \times
\mathbb{C}(q)[[u^{-1}]]\big)^I$ if $v$ satisf\/ies $U_q(\mathbf{L}{\mathfrak g})v = V$, $x_{i,r}^+v = 0$ for all $i \in I$,
$r \in \mathbb{Z}$, and $\phi_i^{\pm}(u) v = \gamma^{\pm}_i(u)v$ for all $i\in I$.
Here $\phi_i^{\pm}(u) \in U_q(\mathbf{L}{\mathfrak h})[[u^{\pm 1}]]$ are def\/ined as follows:
\begin{gather*}
\phi_i^\pm (u)= k_i^{\pm 1} \exp \left(\pm\big(q - q^{-1}\big)\sum\limits_{r = 1}^{\infty} h_{i, \pm r} u^{\pm r}\right).
\end{gather*}

\begin{Theorem}[\cite{MR1357195}] \quad

\begin{enumerate}\itemsep=0pt

\item[{\normalfont(i)}] If $V$ is a~finite-dimensional simple $U_q(\mathbf{L}{\mathfrak g})$-module of type
$1$, then $V$ is highest $\ell$-weight, and its highest $\ell$-weight $\big(\gamma^+_i(u),\gamma_i^-(u)\big)_{i \in I}$
satisfies
\begin{gather}
\label{eq:highest_weight}
\gamma^{\pm}_i(u) = q^{\deg \boldsymbol{\pi}_i(u)}\left(\frac{\boldsymbol{\pi}_i(q^{-1}u)}{\boldsymbol{\pi}_i(qu)}\right)^{\pm}
\end{gather}
for some polynomials $\boldsymbol{\pi}_i(u) \in \mathbb{C}(q)[u]$ whose constant terms are $1$.
Here $(\;)^\pm$ denote the expansions at $u=0$ and $u=\infty$ respectively.

\item[{\normalfont(ii)}] \sloppy Conversely, for every $I$-tuple of polynomials $\boldsymbol{\pi}=\big(\boldsymbol{\pi}_1(u),\ldots,\boldsymbol{\pi}_n(u)\big)$
such that \mbox{$\boldsymbol{\pi}_i(0) =1$}, there exists a~unique $($up to isomorphism$)$ finite-dimensional
simple highest $\ell$-weight $U_q(\mathbf{L}{\mathfrak g})$-module of type $1$ with highest $\ell$-weight
$\big(\gamma^+_i(u), \gamma^-_i(u)\big)_{i \in I}$ satisfying~\eqref{eq:highest_weight}.
\end{enumerate}
\end{Theorem}

The $I$-tuple of polynomials $\boldsymbol{\pi}=\big(\boldsymbol{\pi}_1(u),\ldots,\boldsymbol{\pi}_n(u)\big)$ are called \textit{Drinfeld
polynomials}, and we will say by abuse of terminology that the highest $\ell$-weight of $V$ is $\boldsymbol{\pi}$ if the highest
$\ell$-weight $\big(\gamma^+_i(u), \gamma^-_i(u)\big)_{i \in I}$ of $V$ satisf\/ies~\eqref{eq:highest_weight}.
We denote by $L_q(\boldsymbol{\pi})$ the f\/inite-dimensional simple $U_q(\mathbf{L}{\mathfrak g})$-module of type $1$ with
highest $\ell$-weight $\boldsymbol{\pi}$, and by $v_{\boldsymbol{\pi}}$ a~highest $\ell$-weight vector of $L_q(\boldsymbol{\pi})$.

Let $i \mapsto \bar{i}$ be the bijection $I \to I$ determined by $\alpha_{\bar{i}} = -w_\circ(\alpha_i)$.

\begin{Lemma}
[\cite{MR1300632}]
\label{Lem:dual}
For any Drinfeld polynomials $\boldsymbol{\pi}$ we have
\begin{gather*}
L_q(\boldsymbol{\pi})^* \cong L_q(\boldsymbol{\pi}^*)
\end{gather*}
as $U_q(\mathbf{L}{\mathfrak g})$-modules, where $\boldsymbol{\pi}^* = \big(\boldsymbol{\pi}_{\bar{i}}\big(q^{-h^\vee}u\big)\big)_{i\in
I}$ and $h^\vee$ is the dual Coxeter number.
\end{Lemma}

\subsection{Minimal af\/f\/inizations}\label{subsection:minimal_aff}

For an $I$-tuple of polynomials $\boldsymbol{\pi}=\big(\boldsymbol{\pi}_i(u)\big)_{i \in I}$, set $\mathrm{wt}(\boldsymbol{\pi}) =
\sum\limits_{i\in I}\varpi_i \deg \boldsymbol{\pi}_i \in P^+$.

\begin{Definition}[\cite{MR1367675}]  Let $\lambda \in P^+$.
\begin{enumerate}\itemsep=0pt
\item[(i)] A simple f\/inite-dimensional $U_q(\mathbf{L}{\mathfrak g})$-module $L_q(\boldsymbol{\pi})$ is said to be an
\textit{affinization} of~$V_q(\lambda)$ if~$\mathrm{wt} (\boldsymbol{\pi}) = \lambda$.

\item[(ii)] Af\/f\/inizations $V$ and $W$ of~$V_q(\lambda)$ are said to be \textit{equivalent} if they are isomorphic as
$U_q({\mathfrak g})$-modules.
We denote by~$[V]$ the equivalence class of~$V$.
\end{enumerate}
\end{Definition}

If $V$ is an af\/f\/inization of $V_q(\lambda)$, as a~$U_q({\mathfrak g})$-module we have
\begin{gather*}
V \cong V_q(\lambda) \oplus \bigoplus_{\mu < \lambda} V_q(\mu)^{\oplus m_{\mu}(V)}
\end{gather*}
with some $m_{\mu}(V) \in \mathbb{Z}_{\ge 0}$.
Let $V$ and $W$ be af\/f\/inizations of $V_q(\lambda)$, and def\/ine $m_{\mu}(V)$, $m_{\mu}(W)$ as above.
We write $[V] \le [W]$ if for all $\mu \in P^+$, either of the following holds:
\begin{enumerate}\itemsep=0pt
\item[(i)] $m_\mu(V) \le m_\mu(W)$, or
\item[(ii)] there exists some $\nu > \mu$ such that $m_{\nu}(V) < m_{\nu}(W)$.
\end{enumerate}
Then $\le$ def\/ines a~partial ordering on the set of equivalence classes of af\/f\/inizations of
$V_q(\lambda)$~\cite[Proposition~3.7]{MR1367675}.

\begin{Definition}[\cite{MR1367675}]
We say an af\/f\/inization $V$ of $V_q(\lambda)$ is \textit{minimal} if $[V]$ is minimal in
the set of equivalence classes of af\/f\/inizations of~$V_q(\lambda)$ with respect to this ordering.
\end{Definition}

For $i \in I$, $a \in \mathbb{C}(q)^\times$ and $m \in \mathbb{Z}_{> 0}$, def\/ine an $I$-tuple of polynomials
$\boldsymbol{\pi}_{m,a}^{(i)}$ by
\begin{gather*}
\big(\boldsymbol{\pi}_{m,a}^{(i)}\big)_j(u) =
\begin{cases}
\big(1-aq^{-m+1}u\big)\big(1-aq^{-m+3}u\big) \cdots \big(1-aq^{m-1}u\big), & j = i,
\\
1, & j \neq i.
\end{cases}
\end{gather*}
We set $\boldsymbol{\pi}_{0,a}^{(i)} = (1,1,\ldots,1)$ for every $i \in I$ and $a \in \mathbb{C}(q)^\times$.
The simple modules $L_q(\boldsymbol{\pi}_{m,a}^{(i)})$ are called \textit{Kirillov--Reshetikhin modules}.

Let us recall the classif\/ication of minimal af\/f\/inizations in the regular case of type $D$, which was given
in~\cite{MR1402568}.
(Similar results also hold in type $E$.
See~\cite{MR1347873} for type $ABCFG$, in which minimal af\/f\/inizations are unique up to equivalence.)
For that, we f\/ix several notation.
Set $S = \{1,n-1,n\} \subseteq I$ and def\/ine the subsets $I_s \subseteq I$ $(s \in S)$ by $I_1 =\{1,2,\ldots, n-3\}$,
$I_{n-1} =\{n-1\}$, $I_n =\{n\}$.
Note that $I_s$ is the connected component of the subdiagram $I \setminus \{n-2\}$ containing~$s$, and $I \setminus I_s$
is the maximal type $A$ subdiagram of $I$ not containing~$s$.
For $s \in S$, $\varepsilon \in \{\pm \}$, $\lambda = \sum\limits_{i \in I} \lambda_i \varpi_i \in P^+$ and $a \in
\mathbb{C}(q)^\times$, def\/ine Drinfeld polynomials $\boldsymbol{\pi}_s^{\varepsilon}(\lambda,a)$ as follows:
\begin{list}{$\circ}{\setlength{\leftmargin}{17pt}} \itemsep=0pt

\item[$\circ$] When $s = 1$, set $\boldsymbol{\pi}_1^{\varepsilon}(\lambda,a) =
\prod\limits_{i \in I} \boldsymbol{\pi}_{\lambda_i,a_i}^{(i)}$ $($the product being def\/ined
component-wise$)$ with $a_1 = a$ and
\begin{gather*}
a_i =
\begin{cases}
aq^{\varepsilon(\lambda_1 + 2\sum\limits_{1 < j < i} \lambda_j + \lambda_i + i-1)}, & 2\le i \le n-2,
\vspace{1mm}\\
aq^{\varepsilon(\lambda_1 + 2\sum\limits_{1 < j <n-1}\lambda_j + \lambda_i + n-2)}, & i =n-1, n.
\end{cases}
\end{gather*}
\item[$\circ$] When $s = n-1$ or $n$, set $\boldsymbol{\pi}_s^{\varepsilon}(\lambda,a) = \prod\limits_{i \in I}
\boldsymbol{\pi}_{\lambda_i,a_i}^{(i)}$ with $a_1 = a$ and
\begin{gather*}
a_i =
\begin{cases}
aq^{\varepsilon(\lambda_1 + 2\sum\limits_{1 < j < i} \lambda_j + \lambda_i + i-1)}, & 2\le i \le n-2,
\vspace{1mm}\\
aq^{\varepsilon(\lambda_1 + 2\sum\limits_{1 < j <n-1}\lambda_j + \lambda_i + n-2)}, & i = s,
\vspace{1mm}\\
aq^{\varepsilon(\lambda_1 + 2\sum\limits_{1 < j < n-2}\lambda_j -\lambda_i + n-4)}, & i \in \{n-1,n\},\quad i \neq s.
\end{cases}
\end{gather*}
\end{list}

\begin{Remark}\label{Rem:restriction}
The Drinfeld polynomials $\boldsymbol{\pi}_s^{\varepsilon}(\lambda,a)$ are determined so that they satisfy the
following property: if $r \in S \setminus \{s\}$ and $J = I \setminus I_r$, the simple $U_q(\mathbf{L}{\mathfrak
g}_J)$-module $L_q\big(\boldsymbol{\pi}_s^{\varepsilon}(\lambda,a)_J\big)$ is a~minimal af\/f\/inization of the $U_q({\mathfrak
g}_J)$-module $V_q(\lambda|_{{\mathfrak h}_J})$, where $\boldsymbol{\pi}_s^{\varepsilon}(\lambda,a)_J$ denotes the $J$-tuple
$\big(\boldsymbol{\pi}_s^{\varepsilon}(\lambda,a)_j(u)\big)_{j \in J}$ (see~\cite[Theorem~3.1]{MR1402568}).
\end{Remark}

Def\/ine the \textit{support} of $\lambda = \sum\limits_i \lambda_i \varpi_i\in P^+$ by
\begin{gather*}
\supp (\lambda) = \{ i \in I \,|\, \lambda_i > 0 \} \subseteq I.
\end{gather*}

\begin{Theorem}[\protect{\cite[Theorem~6.1]{MR1402568}}]
\label{Thm:Classification}
Let $\lambda \in P^+$.
\begin{enumerate}\itemsep=0pt
\item[{\normalfont(i)}] If $\supp(\lambda) \cap I_s = \varnothing$ for some $s \in S$, then there exists a~unique
equivalence class of minimal affinizations of $V_q(\lambda)$, and the equivalence class is given by
\begin{gather*}
\big\{ L_q\big(\boldsymbol{\pi}_r^{\varepsilon}(\lambda,a)\big) \,\big|\, \varepsilon \in \{ \pm\},\, a \in
\mathbb{C}(q)^\times\big\}
\end{gather*}
with $r \in S \setminus \{s\}$ $($here the choice of $r$ is irrelevant since $\boldsymbol{\pi}_r^{\pm}(\lambda,a) =
\boldsymbol{\pi}_{r'}^{\pm}(\lambda,a)$ holds for $r$, $r'\in S\setminus\{s\})$.

\item[{\normalfont(ii)}] If $\supp(\lambda) \cap I_s \neq \varnothing$ for all $s \in S$ and $\lambda_{n-2} > 0$, then
there exist exactly three equivalence classes of minimal affinizations of $V_q(\lambda)$, and for each $s \in S$
\begin{gather*}
\big\{ L_q\big(\boldsymbol{\pi}_s^{\varepsilon}(\lambda,a)\big)\, \big|\,\varepsilon \in \{ \pm\}, \, a \in
\mathbb{C}(q)^\times\big\}
\end{gather*}
forms an equivalence class.
\end{enumerate}
\end{Theorem}

We call $\lambda\in P^+$ \textit{regular} if $\lambda$ satisf\/ies one of the assumptions of (i) or (ii) in
Theorem~\ref{Thm:Classification}.
We call a~minimal af\/f\/inization is regular if its highest weight is regular.

\begin{Remark}[\cite{MR1376937}]  In the remaining case when $\supp(\lambda) \cap I_s \neq \varnothing$ for all $s
\in S$ and $\lambda_{n-2} = 0$, the number of equivalence classes of minimal af\/f\/inizations increases unboundedly with~$\lambda$, and the classif\/ication of minimal af\/f\/inizations has not been given except for the type~$D_4$.
\end{Remark}

\subsection{Classical limits and graded limits}

Let ${\bf A} = \mathbb{C}[q,q^{-1}]$ be the ring of Laurent polynomials with complex coef\/f\/icients, and denote by~$U_{\bf
A}(\mathbf{L}{\mathfrak g})$ the ${\bf A}$-subalgebra of $U_q(\mathbf{L}{\mathfrak g})$ generated by $\big\{ k_i^{\pm
1}, \big(x_{i,r}^{\pm}\big)^{k}/[k]_q! \,|\, i \in I, r \in \mathbb{Z}, k \in \mathbb{Z}_{> 0} \big\}$, where we set
$[k]_q = (q^k-q^{-k})/(q-q^{-1})$ and $[k]_q! = [k]_q [k-1]_q\cdots [1]_q$.
Def\/ine $U_{\mathbf{A}}({\mathfrak g}) \subseteq U_q({\mathfrak g})$ in a~similar way.
We def\/ine $\mathbb{C}$-algebras $U_1(\mathbf{L}{\mathfrak g})$ and $U_1({\mathfrak g})$ by
\begin{gather*}
U_1(\mathbf{L}{\mathfrak g}) =\mathbb{C} \otimes_{\bf A} U_{\bf A}(\mathbf{L}{\mathfrak g})
\qquad
\text{and}
\qquad
U_1({\mathfrak g}) =\mathbb{C} \otimes_{\bf A} U_{\bf A} ({\mathfrak g}),
\end{gather*}
where we identify $\mathbb{C}$ with ${\bf A}/\langle q -1\rangle$.
Then the following $\mathbb{C}$-algebra isomorphisms are known to hold~\cite{MR1227098},~\cite[Proposition
9.3.10]{MR1300632}:
\begin{gather}
\label{eq:isom}
U(\mathbf{L}{\mathfrak g}) \cong U_1(\mathbf{L}{\mathfrak g})/\langle k_i-1\,|\, i \in I\rangle_{U_1(\mathbf{L}{\mathfrak
g})},
\qquad
U({\mathfrak g}) \cong U_1({\mathfrak g}) /\langle k_i-1\,|\, i \in I\rangle_{U_1({\mathfrak g})},
\end{gather}
where $\langle k_i-1\,|\, i \in I\rangle_{U_1(\mathbf{L}{\mathfrak g})}$ denotes the two-sided ideal of
$U_1(\mathbf{L}{\mathfrak g})$ generated by $\{k_i -1 \,|\, i \in I\}$, and $\langle k_i-1\,|\, i \in
I\rangle_{U_1({\mathfrak g})}$ is def\/ined similarly.

Let $\boldsymbol{\pi}=\big(\boldsymbol{\pi}_1(u),\ldots,\boldsymbol{\pi}_n(u)\big)$ be Drinfeld polynomials, and assume that there exists $b \in
\mathbb{C}^\times$ such that all the roots of $\boldsymbol{\pi}_i(u)$'s are contained in the set $b q^\mathbb{Z}$ (it is known
that in order to describe the category of f\/inite-dimensional $U_q(\mathbf{L}{\mathfrak g})$-modules, it is essentially
enough to consider representations attached to such families of Drinfeld polynomials.
For example, see~\cite[Section~3.7]{HL}).
Note that $\boldsymbol{\pi}_s^{\pm}(\lambda,a)$ satisf\/ies this assumptions when $a \in \mathbb{C}^\times q^\mathbb{Z}$.
Let $L_{\bf A}(\boldsymbol{\pi})$ be the $U_{\bf A}(\mathbf{L}{\mathfrak g})$-submodule of $L_q(\boldsymbol{\pi})$ generated by a~highest $\ell$-weight vector $v_{\boldsymbol{\pi}}$.
Then by the isomorphism~\eqref{eq:isom},
\begin{gather*}
\overline{L_q(\boldsymbol{\pi})} = \mathbb{C} \otimes_{\bf A} L_{\bf A}(\boldsymbol{\pi})
\end{gather*}
becomes a~f\/inite-dimensional $\mathbf{L}{\mathfrak g}$-module, which is called the \textit{classical limit} of
$L_q(\boldsymbol{\pi})$.

Def\/ine a~Lie algebra automorphism $\varphi_b\colon {\mathfrak g}[t] \to {\mathfrak g}[t]$ by
\begin{gather*}
\varphi_b\big(x \otimes f(t)\big) = x \otimes f(t - b)
\qquad
\text{for}
\quad
x \in {\mathfrak g}, f \in \mathbb{C}[t].
\end{gather*}
We consider $\overline{L_q(\boldsymbol{\pi})}$ as a~${\mathfrak g}[t]$-module by restriction, and def\/ine a~${\mathfrak
g}[t]$-module $L(\boldsymbol{\pi})$ by the pull-back $\varphi_b^*\big(\overline{L_q(\boldsymbol{\pi})}\big)$.
We call $L(\boldsymbol{\pi})$ the \textit{graded limit} of $L_q(\boldsymbol{\pi})$.
In fact, at least when $\boldsymbol{\pi} = \boldsymbol{\pi}_s^{\pm}(\lambda,a)$, it turns out later from our main theorems that
$L(\boldsymbol{\pi})$ is a~graded ${\mathfrak g}[t]$-module, which justif\/ies the name ``graded limit".
Set $\bar{v}_{\boldsymbol{\pi}} = 1 \otimes_{\mathbf{A}} v_{\boldsymbol{\pi}} \in L(\boldsymbol{\pi})$, which generates $L(\boldsymbol{\pi})$ as
a~${\mathfrak g}[t]$-module.
The following properties of graded limits are easily proved from the construction (see~\cite{MR1836791}).

\begin{Lemma}
\label{Lem:elementary_minimal}
Assume $\mathrm{wt}(\boldsymbol{\pi}) = \lambda$.
\begin{enumerate}\itemsep=0pt
\item[{\normalfont(i)}] There exists a~surjective ${\mathfrak g}[t]$-module homomorphism from $L(\boldsymbol{\pi})$ to $V(\lambda,0)$
mapping $\bar{v}_{\boldsymbol{\pi}}$ to a~highest weight vector.

\item[{\normalfont(ii)}] The vector $\bar{v}_{\boldsymbol{\pi}}$ satisfies the relations
\begin{gather*}
{\mathfrak n}_+[t]\bar{v}_{\boldsymbol{\pi}} = 0,
\qquad
\big(h\otimes t^k\big)\bar{v}_{\boldsymbol{\pi}} = \delta_{k0}\langle h, \lambda \rangle \bar{v}_{\boldsymbol{\pi}}
\qquad
\text{for}
\quad
h \in {\mathfrak h}, k \ge 0,
\qquad
\text{and}
\\
f_i^{\langle\alpha_i^\vee, \lambda\rangle+1}\bar{v}_{\boldsymbol{\pi}} = 0
\qquad
\text{for}
\quad
i \in I.
\end{gather*}

\item[{\normalfont(iii)}] We have
\begin{gather*}
\ch L_q(\boldsymbol{\pi}) = \ch L(\boldsymbol{\pi}).
\end{gather*}

\item[{\normalfont(iv)}] For every $\mu \in P^+$ we have
\begin{gather*}
\big[L_q(\boldsymbol{\pi}): V_q(\mu)\big] = \big[L(\boldsymbol{\pi}): V(\mu)\big],
\end{gather*}
where the left- and right-hand sides are the multiplicities as a~$U_q({\mathfrak g})$-module and ${\mathfrak g}$-module,
respectively.
\end{enumerate}
\end{Lemma}

\section{Main theorems and corollaries}\label{Section3}

Throughout this section, we f\/ix $\lambda = \sum\limits_{i \in I} \lambda_i \varpi_i \in P^+$, $\varepsilon \in \{\pm\}$
and $a \in \mathbb{C}^\times q^\mathbb{Z}$, and abbreviate $\boldsymbol{\pi}_s = \boldsymbol{\pi}_s^\varepsilon(\lambda,a)$ for $s \in S
= \{1,n-1,n\}$.

\subsection{Main theorems}\label{subsection:main_theorems}

Denote by $\tau_{0,1} \in \Sigma$ (resp.~$\tau_{n-1,n} \in \Sigma$) the diagram automorphism interchanging the nodes $0$
and~$1$ (resp.~$n-1$ and $n$).
We will not use other elements of $\Sigma$ in the sequel.

For $s \in S$ and $1 \le j \le n$, def\/ine $\xi_j^{(s)}=\xi_j^{(s)}(\lambda) \in \widehat{P}$ as follows:
\begin{list}{$\circ$}{\setlength{\leftmargin}{12pt}\setlength{\itemsep}{7pt}}\itemsep=0pt
\item
When $s=1$, let $m, m'$ be such that $\{m,m'\} = \{n-1,n\}$, $\lambda_m = \max\{\lambda_{n-1}, \lambda_n\}$ and
$\lambda_{m'} = \min\{\lambda_{n-1},\lambda_n\}$, and def\/ine
\begin{gather*}
\xi_j^{(1)} =
\begin{cases}
\lambda_j(\varpi_j + \Lambda_0), & 1 \le j \le n-2,
\\
\lambda_{m'}(\varpi_{n-1} + \varpi_n + \Lambda_0), & j=n-1,
\\
(\lambda_m - \lambda_{m'})(\varpi_m + \Lambda_0), & j = n.
\end{cases}
\end{gather*}
\item
When $s=n$, set
\begin{gather*}
\ell =
\begin{cases}
0, & \text{if} \quad \sum\limits_{i=1}^{n-3} \lambda_i < \lambda_{n-1},
\\
\max \Big\{1\le j \le n-3 \,|\, \sum\limits_{i=j}^{n-3}\lambda_i \ge \lambda_{n-1}\Big\}, & \text{otherwise},
\end{cases}
\end{gather*}
and $\bar{\lambda} = \lambda_{n-1} - \sum\limits_{i=\ell+1}^{n-3}\lambda_i$.
Then def\/ine
\begin{gather*}
\xi_j^{(n)} =
\begin{cases}
\lambda_j(\varpi_j + \Lambda_0), & 1 \le j <\ell, \quad j=n-2,n,
\\
\lambda_\ell(\varpi_\ell + \Lambda_0) + \bar{\lambda}\varpi_{n-1}, &  j=\ell,
\\
\lambda_j(\varpi_j + \varpi_{n-1} + \Lambda_0) + \delta_{\ell0}\delta_{j1}\bar{\lambda}(\varpi_{n-1} + \Lambda_0), &
\ell<j<n-2,
\\
0, & j=n-1.
\end{cases}
\end{gather*}
\item
When $s=n-1$, set $\xi_j^{(n-1)}(\lambda) = \tau_{n-1,n}\big(\xi_j^{(n)}(\tau_{n-1,n}\lambda)\big)$.
\end{list}

Note that we have $\lambda \equiv \sum\limits_{1 \le j \le n}\xi_j^{(s)}$ mod $\mathbb{Z} \Lambda_0 + \mathbb{Q}\delta$
for all $s \in S$.

\begin{Theorem}
\label{Thm:Main2}
The graded limit $L(\boldsymbol{\pi}_s)$ is isomorphic to $D\big(w_\circ\xi_1^{(s)},\ldots,w_\circ \xi_{n}^{(s)}\big)$ as
a~${\mathfrak g}[t]$-module.
\end{Theorem}

For $\alpha = \sum\limits_{i \in I} n_i \alpha_i \in \Delta_+$, set $\supp(\alpha) = \{ i \in I\,|\, n_i > 0\}
\subseteq I$.
We def\/ine a~subset $\Delta^{(s)}_+ \subseteq \Delta_+$ for $s \in S$ by
\begin{gather*}
\Delta^{(s)}_+ = \bigcup_{r \in S \setminus \{s\}} \big\{\alpha \in \Delta_+ \,|\, \supp(\alpha) \subseteq I
\setminus I_r \big\}.
\end{gather*}
Note that if $\alpha \in \Delta^{(s)}_+$, then the coef\/f\/icient of $\alpha_i$ in $\alpha$ is $0$ or $1$ for all $i \in
I$.

\begin{Theorem}
\label{Thm:Theorem1}
The graded limit $L(\boldsymbol{\pi}_s)$ is isomorphic to the cyclic ${\mathfrak g}[t]$-module generated by a~nonzero vector $v$
subject to the relations
\begin{gather*}
{\mathfrak n}_+[t]v = 0,
\qquad
\big(h \otimes t^k\big)v= \delta_{k0}\langle h, \lambda \rangle v
\qquad
\text{for}
\quad
h \in {\mathfrak h}, \quad k \ge 0,
\\
f_i^{\lambda_i+1}v = 0
\qquad
\text{for}
\quad
i \in I
\qquad
\text{and}
\qquad
(f_{\alpha}\otimes t)v= 0
\qquad
\text{for}
\quad
\alpha \in \Delta_+^{(s)}.
\end{gather*}
\end{Theorem}

We prove Theorems~\ref{Thm:Main2} and~\ref{Thm:Theorem1} in Section~\ref{section:Proof}.

\begin{Remark}
 The def\/ining relations of Theorem~\ref{Thm:Theorem1} were conjectured in~\cite[Section~5.11]{MR2587436},
and proved there for ${\mathfrak g}$ of type $D_4$.
Let $I_s' = I_s \sqcup \{n-2\}$, $\lambda_{I_s'} = \sum\limits_{i \in I_s'} \lambda_i \varpi_i$
and $\lambda_{I\setminus I_s'} = \lambda - \lambda_{I_s'}$.
In \textit{loc.~cit.}, the author also conjectured that the graded limit $L(\boldsymbol{\pi}_s)$ is isomorphic to the  
${\mathfrak g}[t]$-submodule of
\begin{gather*}
L\big(\boldsymbol{\pi}_s^{\varepsilon}(\lambda_{I_s'},a)\big) \otimes L\big(\boldsymbol{\pi}_s^{\varepsilon}(\lambda_{I \setminus
I_s'},a)\big)
\end{gather*}
generated by the tensor product of highest weight vectors.
This is easily deduced from Theorem~\ref{Thm:Main2}.
\end{Remark}

\subsection{Corollaries}

The module $D\big(w_\circ \xi_1^{(s)},\ldots,w_\circ \xi_n^{(s)}\big)$ in Theorem~\ref{Thm:Main2} has another
realization introduced in Section~\ref{subsection:b-module}.
Def\/ine $\sigma \in \widetilde{W}$ by
\begin{gather*}
\sigma = \tau_{0,1}\tau_{n-1,n}s_1s_2\cdots s_{n-1}.
\end{gather*}
The proof of the following lemma is straightforward.

\begin{Lemma}\label{Lem;Weyl}\quad
\begin{enumerate}\itemsep=0pt
\item[{\normalfont(i)}] For $0 \le j \le n$, we have
\begin{gather*}
\sigma (\varpi_j + \Lambda_0) \equiv
\begin{cases}
\varpi_{j+1} + \Lambda_0, & 0 \le j \le n-3,
\\
\varpi_{n-1} + \varpi_n + \Lambda_0, & j=n-2,
\\
\varpi_{n-1} + \varpi_1 + \Lambda_0, & j=n-1,
\\
\varpi_{n-1} +\Lambda_0, & j=n,
\end{cases}
\qquad
\mod \mathbb{Q}\delta,
\end{gather*}
and $\sigma(\varpi_{n-1}) \equiv \varpi_{n-1}$ $\mod\mathbb{Q}\delta$.

\item[{\normalfont(ii)}] We have $\ell(w_\circ \sigma^{n-1}) = \ell(w_\circ) + (n-1)\ell(\sigma)$.
\end{enumerate}
\end{Lemma}

Assume $s \neq n-1$ for a~while.
For $1 \le j \le n-1$ def\/ine $\Lambda^{(s)}_j = \sigma^{-j}\xi_j^{(s)}$, and set $\Lambda^{(s)}_n = \xi_n^{(s)}$.
The following assertions are easily checked using Lemma~\ref{Lem;Weyl}(i):
\begin{gather*}
\Lambda^{(1)}_j \equiv
\begin{cases}
\lambda_j\Lambda_0, &  1 \le j \le n-2,
\\
\lambda_{m'}\Lambda_0, & j =n-1,
\end{cases}
\qquad
\mod
\mathbb{Q} \delta,
\qquad
\text{and}
\\
\Lambda^{(n)}_j \equiv
\begin{cases}
\lambda_j\Lambda_0, & 1 \le j < \ell \quad \text{or} \quad j = n-2,
\\
\lambda_\ell\Lambda_0 + \bar{\lambda}\varpi_{n-1}, &  j=\ell,
\\
\lambda_j(\varpi_{n-1}+\Lambda_0) + \delta_{\ell 0}\delta_{j1} \bar{\lambda}(\varpi_n+\Lambda_0), &  \ell< j < n-2,
\end{cases}
\qquad\!\!\!
\mod \mathbb{Q}\delta.
\end{gather*}
In particular, each $\Lambda_{j}^{(s)}$ belongs to $\widehat{P}^+$.
We obtain the following chain of isomorphisms from Theorem~\ref{Thm:Main2}, Lemma~\ref{Lem;Weyl}(ii), and
Proposition~\ref{Prop:construction}:
\begin{gather}
\nonumber
L(\boldsymbol{\pi}_s) \cong D\big(w_\circ\xi_1^{(s)},\ldots,w_\circ \xi_{n}^{(s)}\big) \cong D\big(w_\circ\xi_n^{(s)},
w_\circ\xi_1^{(s)}, \ldots,w_\circ\xi_{n-1}^{(s)}\big)
\\
\phantom{L(\boldsymbol{\pi}_s)}
\cong F_{w_\circ}\big(D\big(\Lambda_n^{(s)}\big)\otimes F_{\sigma}\big(D\big(\Lambda_1^{(s)}\big) \otimes \dots \otimes
F_{\sigma}\big(D\big(\Lambda_{n-2}^{(s)}\big) \otimes F_{\sigma} D\big(\Lambda_{n-1}^{(s)}\big)\big)\cdots \big)\big),
\label{eq:isomorphism}
\end{gather}
where the second isomorphism obviously holds by def\/inition.
Hence by Proposition~\ref{Prop:character} and Lemma~\ref{Lem:elementary_minimal}(iii), the following holds.

\begin{Corollary}
\label{Cor:characters}
If $s \in \{1,n\}$, we have
\begin{gather*}
\ch L_q(\boldsymbol{\pi}_s) = \mathcal{D}_{w_\circ}\big(e^{\Lambda_{n}^{(s)}}\cdot
\mathcal{D}_{\sigma}\big(e^{\Lambda_1^{(s)}} \cdots \mathcal{D}_{\sigma}\big(e^{\Lambda_{n-2}^{(s)}}\cdot
\mathcal{D}_{\sigma}\big(e^{\Lambda_{n-1}^{(s)}}\big) \big)\cdots\big)\big)\big|_{e^{\Lambda_0} =e^\delta= 1}.
\end{gather*}
\end{Corollary}

Let $\lambda' = \tau_{n-1,n}\lambda$, and set $\boldsymbol{\pi}_{n}' = \boldsymbol{\pi}_{n}^{\varepsilon}(\lambda',a)$.
It is easily seen from Theorem~\ref{Thm:Main2} that
\begin{gather*}
\ch L_q(\boldsymbol{\pi}_{n-1}) = \tau_{n-1,n}\ch L_q(\boldsymbol{\pi}_n').
\end{gather*}
Hence we also obtain the character in the case $s = n-1$.

\begin{Remark}\quad
\begin{enumerate}\itemsep=0pt
\item[(i)] It is possible to use other elements of $\widetilde{W}$ in the expression of $\ch L_q(\boldsymbol{\pi}_s)$.
That is, if $w_j \in \widetilde{W}$ $(1 \le j \le n-1)$ satisfy $w_{[1,j]}\Lambda_j^{(s)} = \xi_j^{(s)}$ and
$\ell(w_\circ w_{[1,n-1]}) = \ell(w_\circ) + \sum\limits_{j=1}^{n-1} \ell(w_j)$ (here we set $w_{[1,j]} = w_1w_2\cdots
w_j$), then it follows that
\begin{gather*}
\ch L_q(\boldsymbol{\pi}_s) = \mathcal{D}_{w_\circ}\big(e^{\Lambda_{n}^{(s)}}\cdot
\mathcal{D}_{w_1}\big(e^{\Lambda_1^{(s)}} \cdots \mathcal{D}_{w_{n-2}}\big(e^{\Lambda_{n-2}^{(s)}}\cdot
\mathcal{D}_{w_{n-1}}\big(e^{\Lambda_{n-1}^{(s)}}\big) \big)\cdots\big)\big)\big|_{e^{\Lambda_0} =e^\delta= 1}.
\end{gather*}
For example $w_j = s_{j-1}s_{j-2}\cdots s_1 \tau_{0,1}$ satisfy the above conditions when $s=1$.
Our choice is made so that the results are stated in a~uniform way.

\item[(ii)] The right-hand side of the isomorphism~\eqref{eq:isomorphism} has a~crystal analog, and using this we can express
the multiplicities of~$L_q(\boldsymbol{\pi}_s)$ in terms of crystal bases.
For the details, see~\cite[Corollary~4.11]{N}.
\end{enumerate}
\end{Remark}

Our next result is a~formula for multiplicities of simple f\/inite-dimensional $U_q({\mathfrak g})$-modules in~$L_q(\boldsymbol{\pi}_1)$ which can be deduced from our Theorem~\ref{Thm:Theorem1} and the results of~\cite{MR2763623}
and~\cite{Sam}.
For that, we prepare a~lemma.

\begin{Lemma}
\label{Lem:multiplicity}
Assume that $V$ is a~cyclic finite-dimensional ${\mathfrak g}[t]$-module generated by a~${\mathfrak h}$-weight vector
$v$, and ${\mathfrak n}_+[t] \oplus t{\mathfrak h}[t]$ acts trivially on $v$.
Let $\mu \in P^+$, and $W$ be the ${\mathfrak g}[t]$-submodule of $V \otimes V(\mu,0)$ generated by $v \otimes v_\mu$,
where $v_\mu$ denotes a~highest weight vector.
Then for every $\nu \in P^+$, we have
\begin{gather*}
\big[W: V(\nu + \mu)\big] = \big[V: V(\nu)\big],
\end{gather*}
where $[\,: \,]$ denotes the multiplicity as a~${\mathfrak g}$-module.
\end{Lemma}

\begin{proof}
Note that
\begin{gather}
\label{eq:multi}
\big[W: V(\nu + \mu)\big] = \dim\{w \in W_{\nu + \mu} \,|\, {\mathfrak n}_+w=0 \}.
\end{gather}
Since
\begin{gather*}
W = U\big({\mathfrak n}_-[t]\big)(v \otimes v_\mu) = U({\mathfrak n}_-)\big(U\big(t{\mathfrak n}_-[t]\big) v\otimes
v_\mu\big)
\end{gather*}
and $W$ is a~f\/inite-dimensional ${\mathfrak g}$-module, we see that
\begin{gather*}
\{w \in W_{\nu + \mu} \,|\, {\mathfrak n}_+w = 0\} = \big\{ w \in \big(U\big(t{\mathfrak n}_-[t]\big)v \otimes
v_\mu\big)_{\nu + \mu} \,\big|\, {\mathfrak n}_+w = 0\big\}
\\
\qquad
=\big\{ w \in \big(U\big(t{\mathfrak n}_-[t]\big)v\big)_\nu \,|\, {\mathfrak n}_+w = 0\big\} \otimes v_\mu
=\big\{ w \in \big(U\big({\mathfrak n}_-[t]\big)v\big)_\nu\,|\, {\mathfrak n}_+w=0\big\} \otimes v_\mu.
\end{gather*}
Hence the assertion follows from~\eqref{eq:multi}.
\end{proof}

Let $\mathfrak{sp}_{2n-2}$ be the simple Lie algebra of type $C_{n-1}$, and denote by $P_{\mathfrak{sp}}$ its weight
lattice and by $\varpi_i^{\mathfrak{sp}}$ $(1 \le i \le n-1)$ its fundamental weights.
We assume that $\varpi_i^{\mathfrak{sp}}$ are labeled as~\cite[Section~4.8]{MR1104219}.
Def\/ine a~map $\iota\colon P^+ \to P_{\mathfrak{sp}}^+$ by
\begin{gather*}
\iota\left(\sum\limits_{1\le i \le n} \mu_i \varpi_i\right) = \sum\limits_{1 \le i \le n-2} \mu_i \varpi_i^{\mathfrak{sp}}
+ \min\{\mu_{n-1}, \mu_n \} \varpi_{n-1}^{\mathfrak{sp}}.
\end{gather*}

\begin{Corollary}
\label{Cor:Sam}
For every $\mu \in P^+$, we have
\begin{gather*}
\big[ L_q(\boldsymbol{\pi}_1): V_q(\mu)\big] =
\begin{cases}
\big[ \mathbf{S}_{\iota(\lambda)}\big(V^{\mathfrak{sp}}(\varpi_1^{\mathfrak{sp}})\big):
V^{\mathfrak{sp}}\big(\iota(\mu)\big)\big], & \text{if} \quad \mu_{n} - \mu_{n-1} = \lambda_n - \lambda_{n-1},
\\
0, & \text{otherwise}.
\end{cases}
\end{gather*}
Here $\mathbf{S}_{\nu}$ $(\nu \in P^+_{\mathfrak{sp}})$ denotes the Schur functor {\rm (see
\cite[Section~1]{Sam})} with respect to the partition $\Big(\sum\limits_{j=1}^{n-1} \nu_j,
\sum\limits_{j=2}^{n-1} \nu_j, \ldots, \nu_{n-1}\Big)$, and $V^{\mathfrak{sp}}(\nu)$ denotes the simple
$\mathfrak{sp}_{2n-2}$-module with highest weight $\nu$.
\end{Corollary}

\begin{proof}
It suf\/f\/ices to show that the right-hand side is equal to $\big[ L(\boldsymbol{\pi}_1)\!:\! V(\mu)\big]$ by
Lemma~\ref{Lem:elementary_minimal}(iv).
Note that Theorem~\ref{Thm:Theorem1} and~\cite[Theorem~1]{MR2763623} imply that the graded limit $L(\boldsymbol{\pi}_1)$ is
isomorphic to the ${\mathfrak g}[t]$-module ``$P(\lambda,0)^{\Gamma(\lambda,\Psi)}$'' in the notation
of~\cite{MR2763623}, where we set \mbox{$\Psi \!=\! \{ \alpha \!\in\! \Delta_+ \:|\: (\alpha, \varpi_{n-1} \!+\! \varpi_n) \!=\!2\}$}.
Then in the case $\lambda_{n-1} = \lambda_n$, our assertion is a~consequence of~\cite[Theorem~2]{MR2763623}
and~\cite[Theo\-rem~1]{Sam}.

Let us assume $\lambda_{n-1} \neq \lambda_{n}$, and set $\lambda' = \lambda_m - \lambda_{m'}$.
We have
\begin{gather*}
L(\boldsymbol{\pi}_1) \cong D\big(w_\circ \xi_1^{(1)},\ldots,w_\circ \xi_n^{(1)}\big)
\end{gather*}
by Theorem~\ref{Thm:Main2}, and we easily see that $D\big(w_\circ \xi_n^{(1)}\big) \cong V(\lambda'\varpi_m, 0)$ holds.
Hence by applying Lemma~\ref{Lem:multiplicity} with $V = D\big(w_\circ \xi_1^{(1)},\ldots, w_\circ \xi_{n-1}^{(1)})$ and
$\mu = \lambda'\varpi_m$, we have for every $\nu \in P^+$ that
\begin{gather*}
\Big[ D\big(w_\circ \xi_1^{(1)},
\dots, w_\circ \xi_n^{(1)}\big): V(\nu + \lambda'\varpi_m)\big] = \big[ D\big(w_\circ
\xi_1^{(1)},\dots, w_\circ \xi_{n-1}^{(1)}\big): V(\nu)\big],
\end{gather*}
and the right-hand side is equal to $[L\big(\boldsymbol{\pi}_1^{\varepsilon}(\lambda - \lambda'\varpi_m,a)\big):V(\nu)]$ by Theorem~\ref{Thm:Main2}.
Hence the assertion is deduced from the case $\lambda_{n-1} = \lambda_n$.
The proof is complete.
\end{proof}

\section{Proofs of main theorems}\label{section:Proof}

Note that the theorems for $s=n-1$ and $s=n$ are equivalent because of the existence of diagram automorphism of
${\mathfrak g}$ interchanging $n-1$ and $n$.
Therefore, throughout this section we assume that $s \neq n-1$, and prove the theorems only for the case $s = 1, n$.
Similarly as the previous section we f\/ix $\lambda = \sum\limits_{i \in I} \lambda_i \varpi_i \in P^+$, $\varepsilon \in
\{\pm\}$ and $a \in \mathbb{C}^\times q^\mathbb{Z}$, and write $\boldsymbol{\pi}_s = \boldsymbol{\pi}_s^{\varepsilon}(\lambda,a)$.
Let $M_s(\lambda)$ denote the module def\/ined in Theorem~\ref{Thm:Theorem1}.
We shall prove the existence of three surjective ${\mathfrak g}[t]$-module homomorphisms.
More presicely, we prove $M_s(\lambda) \twoheadrightarrow L(\boldsymbol{\pi}_s)$ in Section~\ref{subsection:3-1},
$D\big(w_\circ \xi^{(s)}_1,\ldots w_\circ \xi^{(s)}_{n}\big) \twoheadrightarrow M_s(\lambda)$ in
Section~\ref{subsection:3-2}, and $L(\boldsymbol{\pi}_s) \twoheadrightarrow D\big(w_\circ \xi^{(s)}_1,\ldots w_\circ\xi^{(s)}_{n}\big)$
in Section~\ref{subsection:3-3}.
Then both Theorems~\ref{Thm:Main2} and~\ref{Thm:Theorem1} immediately follow from them.

\subsection[Proof for $M_s(\lambda)\twoheadrightarrow L(\pi_s)$]{Proof for $\boldsymbol{M_s(\lambda)\twoheadrightarrow L(\pi_s)}$}
\label{subsection:3-1}

Though the proof is similar to that in~\cite{MR2587436}, we will give it for completeness.

Let $v_{\boldsymbol{\pi}_s}$ be a~highest $\ell$-weight vector of $L_q(\boldsymbol{\pi}_s)$, and set $\bar{v}_{\boldsymbol{\pi}_s}=1 \otimes
v_{\boldsymbol{\pi}_s} \in L(\boldsymbol{\pi}_s)$.
In order to prove $M_s(\lambda) \twoheadrightarrow L(\boldsymbol{\pi}_s)$, it is enough to show the relations{\samepage
\begin{gather}
\label{eq:relations}
(f_\alpha \otimes t)\bar{v}_{\boldsymbol{\pi}_s} = 0
\qquad
\text{for}
\quad
\alpha \in \Delta_+^{(s)},
\end{gather}
since the other relations hold by Lemma~\ref{Lem:elementary_minimal}(ii).}

Let $r \in S \setminus \{s\}$ and $J = I \setminus I_r$.
Then the subalgebra $U_q(\mathbf{L}{\mathfrak g}_{J}) \subseteq U_q(\mathbf{L}{\mathfrak g})$ is a~quantum loop algebra
of type $A$.
By {\cite[Lemma 2.3]{MR1402568}}, the $U_q(\mathbf{L}{\mathfrak g}_J)$-submodule of $L_q(\boldsymbol{\pi}_s)$ generated by
$v_{\boldsymbol{\pi}_s}$ is isomorphic to the simple $U_q(\mathbf{L}{\mathfrak g}_{J})$-module with highest $\ell$-weight
$\big((\boldsymbol{\pi}_s)_i(u)\big)_{i \in J}$.
Denote this $U_q(\mathbf{L}{\mathfrak g}_J)$-submodule by $L_q'$.
Then we see from Remark~\ref{Rem:restriction} and~\cite[Theorem 3.1]{MR1402568} that $L_q'$ is also simple as
a~$U_q({\mathfrak g}_{J})$-module.
From this and the construction of graded limits, it follows that the $\mathbf{L}{\mathfrak g}_{J}$-submodule
\begin{gather*}
L'=U(\mathbf{L}{\mathfrak g}_{J})\bar{v}_{\boldsymbol{\pi}_s} \subseteq L(\boldsymbol{\pi}_s)
\end{gather*}
is simple as a~${\mathfrak g}_{J}$-module.
Hence the restriction of the surjective homomorphism $L(\boldsymbol{\pi}_s) \twoheadrightarrow V(\lambda,0)$ in
Lemma~\ref{Lem:elementary_minimal}(i) to $L'$ is an isomorphism, which obviously implies $({\mathfrak g}_J \otimes
t)\bar{v}_{\boldsymbol{\pi}_s}= 0$.
Now the relations~\eqref{eq:relations} obviously follow from the def\/inition of $\Delta_+^{(s)}$.
The proof is complete.

\subsection[Proof for $D\big(w_\circ \xi_1^{(s)},\ldots,w_\circ \xi_{n}^{(s)}\big) \twoheadrightarrow M_s(\lambda)$]{Proof
for $\boldsymbol{D\big(w_\circ \xi_1^{(s)},\ldots,w_\circ \xi_{n}^{(s)}\big) \twoheadrightarrow M_s(\lambda)}$}
\label{subsection:3-2}

Throughout this subsection, we assume that $s \in \{1,n\}$ is f\/ixed.
Note that some notation appearing below may depend on $s$ though it is not written explicitly.

Let us prepare several notation.
For $1 \le p \le q \le n$, set
\begin{gather*}
\alpha_{p,q} =
\begin{cases}
\alpha_p + \alpha_{p+1} + \cdots + \alpha_q, & q \le n-1,
\\
\alpha_p + \alpha_{p+1} + \cdots + \alpha_{n-2} + \alpha_n, & q = n.
\end{cases}
\end{gather*}
Note that
\begin{gather*}
\Delta_+ = \{\alpha_{p,q} \,|\, p\le q, (p,q) \neq (n-1,n)\}\sqcup \{ \alpha_{p,n} + \alpha_{q,n-1}\,|\, p <q < n\}.
\end{gather*}

Set $\sigma_i = s_{i} s_{i+1} \cdots s_{n-1} \in \widehat{W}$ for $1 \le i \le n$ and $\sigma_0 =
\tau_{0,1}\tau_{n-1,n}\sigma_1=\sigma$.
For $0 \le i \le n$ and $1 \le j \le n-1$, def\/ine $\rho_{i,j}\colon \widehat{\Delta} \to \mathbb{Z}_{\ge 0}$ by
\begin{gather*}
\rho_{i,j}(\alpha) = \sum\limits_{k =j}^{n-1} \max \big\{0, -\big( \alpha,
\sigma_i\sigma^{k-j}\Lambda_k^{(s)}\big)\big\}.
\end{gather*}
When $j<n-1$, we have
\begin{gather}
\label{eq:ind}
\rho_{n,j}(\alpha) = \rho_{0,j+1}(\alpha) + \max\big\{0, -\big( \alpha,\Lambda^{(s)}_j\big)\big\} = \rho_{0,j+1}(\alpha)
\qquad
\text{for} \quad \alpha \in \widehat{\Delta}^{\mathrm{re}}_+
\end{gather}
since $\Lambda_j^{(s)} \in \widehat{P}^+$.

\begin{Lemma}
\label{Lem:values}
Let $1 \le i \le n$ and $1 \le j \le n-1$, and assume that $\alpha = \beta + k \delta \in
\widehat{\Delta}^{\mathrm{re}}_+$ satisfies $\rho_{i,j}(\alpha) > 0$.
\begin{enumerate}\itemsep=0pt
\item[{\normalfont(i)}] If $i = n$, we have
\begin{gather*}
\beta \in \{ -\alpha_{p,n-1}\,|\, p < n \} \sqcup \big\{ {-}(\alpha_{p,n} + \alpha_{q,n-1})\,|\, p < q <n\big\}.
\end{gather*}
\item[{\normalfont(ii)}] If $1 \le i \le n-1$, we have
\begin{gather*}
\beta \in \{ \alpha_{i,q} \,|\, i \le q <n\}\sqcup \{-\alpha_{p,i-1}\,|\, p<i \} \sqcup \{- \alpha_{p,n} \,|\, p\neq i\}
\\
\qquad
\sqcup \big\{- (\alpha_{p,n} + \alpha_{q,n-1}) \,|\, p <q <n,p \neq i, q \neq i \big\}.
\end{gather*}
\end{enumerate}
\end{Lemma}

\begin{proof}
In both the cases $s=1$ and $s=n$, it follows from Lemma~\ref{Lem;Weyl}(i) that
\begin{gather*}
\sigma^{k-j}\Lambda_k^{(s)} \in \sum\limits_{0 \le p \le n} \mathbb{Z}_{\ge 0} (\varpi_p + \Lambda_0) + \sum\limits_{1
\le p \le n} \mathbb{Z}_{\ge 0}(\varpi_p + \varpi_{n-1} + \Lambda_0)
\end{gather*}
holds for every $1 \le j \le k \le n-1$.
Hence $\rho_{n,j}(\alpha) > 0$ implies $\beta \in -\Delta_+$ and $k \ge 1$, and if $\beta = - \sum\limits_{i \in I}
t_i\alpha_i$ then we have
\begin{gather*}
t_p + t_{n-1} > k \ge 1
\qquad
\text{for some}
\quad
1\le p \le n.
\end{gather*}
This immediately implies the assertion (i).
Note that $\rho_{i,j}(\alpha) = \rho_{i+1,j}(s_i\alpha)$ holds for $1\le i \le n-1$ by the def\/inition of $\rho_{i,j}$.
Hence $\rho_{i,j}(\alpha) > 0$ implies that we have either $\alpha=\alpha_i$ or $\alpha = s_i \gamma$ for some $\gamma
\in \widehat{\Delta}^{\mathrm{re}}_+$ such that $\rho_{i+1,j}(\gamma) > 0$.
From this, the assertion (ii) is easily proved by the descending induction on $i$.
\end{proof}

For $0 \le i \le n$ and $1 \le j \le n-1$, set
\begin{gather*}
D(i,j)= D\big(\sigma_i\Lambda_j^{(s)}, \sigma_i\sigma\Lambda_{j+1}^{(s)},
\ldots,\sigma_i\sigma^{n-j-1}\Lambda^{(s)}_{n-1}\big),
\qquad
\text{and}
\\
v(i,j)= v_{\sigma_i\Lambda_j^{(s)}} \otimes v_{\sigma_i \sigma\Lambda^{(s)}_{j+1}} \otimes \dots \otimes
v_{\sigma_i\sigma^{n-j-1} \Lambda^{(s)}_{n-1}} \in D(i,j).
\end{gather*}
For $\alpha = \beta + k\delta \in \widehat{\Delta}^{\mathrm{re}}$ with $\beta \in \Delta$ and $k \in \mathbb{Z}$, denote
by $x_\alpha \in \widehat{{\mathfrak g}}$ the vector $e_\beta \otimes t^k$.
For $i \in \widehat{I}$, def\/ine a~Lie subalgebra $\widehat{{\mathfrak n}}_i$ of $\widehat{{\mathfrak n}}_+$ by
\begin{gather*}
\widehat{{\mathfrak n}}_i = \bigoplus_{\alpha \in \widehat\Delta_+^{\mathrm{re}} \setminus \{\alpha_i\}} \mathbb{C} x_\alpha
\oplus t{\mathfrak h}[t].
\end{gather*}
We shall determine the generators of the annihilators $\mathrm{Ann}_{U(\widehat{{\mathfrak n}}_+)}v(i,j)$ inductively, along
the lines of~\cite[Section 5.1]{N}.
For that, we need the following lemma which is proved in~\cite[Section 3]{MR826100} (see also~\cite[Lemma 5.3]{N}).

\begin{Lemma}
\label{Lem:proceed}
Let $V$ be an integrable $\widehat{{\mathfrak g}}$-module, $T$ a~f\/inite-dimensional $\widehat{{\mathfrak b}}$-submodule of $V$,
$i \in \widehat{I}$ and $\xi \in \widehat{P}$ such that $( \alpha_i, \xi ) \ge 0$.
Assume that the following conditions hold:
\begin{enumerate}\itemsep=0pt
\item[{\normalfont(i)}] $T$ is generated by a~$\widehat{{\mathfrak h}}$-weight vector $v \in T_\xi$ satisfying $e_i v=0$.

\item[{\normalfont(ii)}] There is an $\ad(e_i)$-invariant left $U(\widehat{{\mathfrak n}}_i)$-ideal $\mathcal{I}$ such that
\begin{gather*}
\mathrm{Ann}_{U(\widehat{{\mathfrak n}}_+)}v = U(\widehat{{\mathfrak n}}_+)e_i + U(\widehat{{\mathfrak n}}_+)\mathcal{I}.
\end{gather*}
\item[{\normalfont(iii)}] We have $\ch_{\widehat{{\mathfrak h}}} F_i T = \mathcal{D}_i \ch_{\widehat{{\mathfrak h}}}
T$.
\end{enumerate}

Let $v' = f_i^{( \alpha_i, \xi )} v$.
Then we have
\begin{gather*}
\mathrm{Ann}_{U(\widehat{{\mathfrak n}}_+)} v'= U(\widehat{{\mathfrak n}}_+)e_i^{(\alpha_i, \xi)+1} + U(\widehat{{\mathfrak
n}}_+)r_i( \mathcal{I}),
\end{gather*}
where $r_i$ denotes the algebra automorphism of $U(\widehat{{\mathfrak g}})$ corresponding to the reflection $s_i$.
\end{Lemma}

\begin{Proposition}
\label{Prop:annihilators}
The following assertion $(\mathrm{A}_{i,j})$ holds for every $0\le i \le n$ and $1 \le j \le n-1$
\begin{gather*}
 (\mathrm{A}_{i,j})
\quad
\mathrm{Ann}_{U(\widehat{{\mathfrak n}}_+)}v(i,j) = U(\widehat{{\mathfrak n}}_+)\left(\sum\limits_{\alpha \in
\widehat{\Delta}^{\mathrm{re}}_+} \mathbb{C} x_{\alpha}^{\rho_{i,j}(\alpha)+1} + t{\mathfrak h}[t]\right).
\end{gather*}
\end{Proposition}

\begin{proof}
We will prove the assertion by the descending induction on $(i,j)$.
The assertion $(\mathrm{A}_{n,n-1})$ is obvious since $D(n,n-1) = \mathbb{C} v(n,n-1)$ is a~trivial $\widehat{{\mathfrak
n}}_+$-module and $\rho_{n,n-1}(\alpha) = 0$ for all $\alpha \in \widehat{\Delta}^{\mathrm{re}}_+$.
Since
\begin{gather*}
v(n,j) = v_{\Lambda_j^{(s)}} \otimes v(0,j+1),
\end{gather*}
we easily see that $(\mathrm{A}_{0,j+1})$ implies $(\mathrm{A}_{n,j})$ by $\eqref{eq:ind}$, and $(\mathrm{A}_{1,j})$
implies $(\mathrm{A}_{0,j})$ since
\begin{gather*}
D(0,j) \cong F_{\tau_{0,1}\tau_{n-1,n}}D(1,j)
\qquad
\text{and}
\qquad
\rho_{0,j}(\alpha) = \rho_{1,j}(\tau_{0,1}\tau_{n-1,n}\alpha).
\end{gather*}
It remains to show that $(\mathrm{A}_{i,j})$ implies $(\mathrm{A}_{i-1,j})$ when $2 \le i \le n$.
Let $\xi(i,j) = \sum\limits_{k=j}^{n-1} \sigma_{i}\sigma^{k-j}\Lambda^{(s)}_k \in \widehat{P}$, which is the weight of
$v(i,j)$.
Since $\big(\alpha_{i-1}, \sigma_{i} \sigma^{k-j}\Lambda_k^{(s)}\big) \ge 0$ holds for all $k\ge j$ by
Lemma~\ref{Lem;Weyl}(ii), we have
\begin{gather*}
e_{i-1} v(i,j) = 0
\qquad
\text{and}
\qquad
f_{i-1}^{(\alpha_{i-1}, \xi(i,j))}v(i,j) \in \mathbb{C}^\times v(i-1,j).
\end{gather*}
In addition, we have $\rho_{i-1,j}(\alpha) = \rho_{i,j}(s_{i-1}\alpha)$ for $\alpha \in
\widehat{\Delta}^{\mathrm{re}}_+$ and in particular
\begin{gather*}
\rho_{i-1,j}(\alpha_{i-1}) = \rho_{i,j}(-\alpha_{i-1})=\big(\alpha_{i-1}, \xi(i,j)\big).
\end{gather*}
Therefore, it suf\/f\/ices to show the $\ad(e_{i-1})$-invariance of the left $U(\widehat{{\mathfrak n}}_{i-1})$-ideal
\begin{gather*}
\mathcal{I}_{i,j} = U(\widehat{{\mathfrak n}}_{i-1})\left(\sum\limits_{\alpha \in \widehat{\Delta}_+^{\mathrm{re}} \setminus
\{ \alpha_{i-1} \}} \mathbb{C} x_\alpha^{\rho_{i,j}(\alpha)+1} + t{\mathfrak h}[t]\right)
\end{gather*}
by Lemma~\ref{Lem:proceed} (note that the condition (iii) holds by Proposition~\ref{Prop:character}).
Since $\rho_{i,j}(\alpha_{i-1} + \mathbb{Z}_{>0}\delta)=0$ holds by Lemma~\ref{Lem:values}, we have
\begin{gather*}
\big[e_{i-1}, t{\mathfrak h}[t]\big] = e_{i-1} \otimes t\mathbb{C}[t] \subseteq \mathcal{I}_{i,j}.
\end{gather*}
Hence it is enough to verify that
\begin{gather}
\label{claim}
\big[e_{i-1}, x_\alpha^{\rho_{i,j}(\alpha)+1}\big] \in \mathcal{I}_{i,j}
\end{gather}
for every $\alpha \in \widehat{\Delta}_+^{\mathrm{re}} \setminus \{ \alpha_{i-1} \}$.
If $\alpha = -\alpha_{i-1} + k\delta$ ($k >0$), then~\eqref{claim} follows from $t{\mathfrak h}[t] \oplus e_{i-1}
\otimes t\mathbb{C}[t] \subseteq \mathcal{I}_{i,j}$.
Hence we may assume that $\alpha$ satisf\/ies $\big[\big[e_{i-1},x_\alpha\big],x_\alpha\big] = 0$.
If $\big[e_{i-1},x_\alpha\big]=0$,~\eqref{claim} is obvious, and otherwise we have
\begin{gather*}
\big[e_{i-1},x_\alpha^{\rho_{i,j}(\alpha)+1}\big] \in \mathbb{C} x_\alpha^{\rho_{i,j}(\alpha)}x_{\alpha + \alpha_{i-1}}.
\end{gather*}
It is directly checked from Lemma~\ref{Lem:values} that if $\beta \in \widehat{\Delta}^{\mathrm{re}}_+$
satisf\/ies $\beta-\alpha_{i-1} \in \widehat{\Delta}^{\mathrm{re}}_+$, then $\rho_{i,j}(\beta) = 0$.
Hence we have $\rho_{i,j}(\alpha + \alpha_{i-1}) = 0$, and~\eqref{claim} follows.
The proof is complete.
\end{proof}

In the sequel we write $\rho = \rho_{0,1}$ for brevity.
Note that we have
\begin{gather*}
\rho(\alpha) =\sum\limits_{1 \le j \le n-1} \max\big\{0, -\big(\alpha, \xi_j^{(s)}\big)\big\}.
\end{gather*}
The following assertions are proved from the def\/inition of $\xi_j^{(s)}$'s by a~direct calculation.
\begin{enumerate}\itemsep=0pt
\item[\normalfont(i)] Assume that $s=1$.
\begin{enumerate}\itemsep=0pt
\item[{\normalfont(a)}] For $\beta + k \delta \in \widehat{\Delta}^{\mathrm{re}}_+$,
$\rho(\beta+k\delta) = 0$ holds unless
\begin{gather*}
-\beta \in \big\{\alpha_{p,n} + \alpha_{q,n-1} \,|\, p < q < n\big\}
\qquad
\text{and}
\qquad
k = 1.
\end{gather*}
\item[{\normalfont(b)}] We have
\begin{gather*}
\rho\big(-(\alpha_{p,n} + \alpha_{q,n-1})+ \delta\big) = \sum\limits_{j = q}^{n-2} \lambda_j + \lambda_{m'}.
\end{gather*}
\end{enumerate}
\item[{\normalfont(ii)}] Assume that $s=n$.
\begin{enumerate}\itemsep=0pt
 \item[{\normalfont(a)}] For $\beta + k \delta \in \widehat{\Delta}^{\mathrm{re}}_+$,
$\rho(\beta+k\delta) = 0$ holds unless
\begin{gather*}
-\beta\in \{\alpha_{p,n-1} \,|\, p \le n-3\}
\qquad
\text{and}
\qquad
k=1,
\qquad
\text{or}
\\
-\beta\in \{\alpha_{p,n} + \alpha_{q,n-1} \,|\, p < q < n\}
\qquad
\text{and}
\qquad
k=1,2.
\end{gather*}
\item[{\normalfont(b)}] We have
\begin{gather*}
\rho(-\alpha_{p,n-1} + \delta) = \min\left\{\sum\limits_{j=p}^{n-3} \lambda_j, \lambda_{n-1}\right\},
\qquad
\text{and}
\\
\rho\big(-(\alpha_{p,n} + \alpha_{q,n-1}) + k\delta\big) =
\begin{cases}
\min\left\{\sum\limits_{j=p}^{n-3} \lambda_j, \lambda_{n-1}\right\} + \sum\limits_{j=q}^{n-2} \lambda_j, & \text{if} \quad k=1,
\\[14pt]%
\min\left\{\sum\limits_{j=q}^{n-3} \lambda_j, \lambda_{n-1}\right\}, & \text{if} \quad k=2.
\end{cases}
\end{gather*}
\end{enumerate}
\end{enumerate}
Set $D = D\big(\xi_n^{(s)}, \xi_1^{(s)}, \ldots, \xi_{n-1}^{(s)}\big)$ and $v_D = v_{\xi_n^{(s)}} \otimes v(0,1) \in D$.
By Proposition~\ref{Prop:annihilators}, we have
\begin{gather*}
\mathrm{Ann}_{U(\widehat{{\mathfrak n}}_+)} v_D = \mathrm{Ann}_{U(\widehat{{\mathfrak n}}_+)} v(0,1) = U(\widehat{{\mathfrak
n}}_+)\left(\sum\limits_{\alpha \in \widehat{\Delta}^{\mathrm{re}}_+} \mathbb{C} x_\alpha^{\rho(\alpha)+1} + t{\mathfrak
h}[t]\right).
\end{gather*}
Let $v_M \in M_s(\lambda)$ denote the generator in the def\/inition.
The proof of the following lemma is elementary (see, e.g., \cite[Lemma~4.5]{MR2855081}).

\begin{Lemma}
\label{Lem:el}
If $\beta, \gamma \in \widehat{\Delta}^{\mathrm{re}}_+$ satisfy $(\beta,\gamma) = -1$, $x_\beta^{b+1}v_M =0$ and
$x_\gamma^{c+1}v_M=0$ with $b,c \in \mathbb{Z}_{\ge 0}$, then $x_{\beta + \gamma}^{b+c+1}v_M = 0$ holds.
\end{Lemma}

\begin{Lemma}
\label{Lem:hom}
There exists a~$({\mathfrak h} \oplus \widehat{{\mathfrak n}}_+)$-module homomorphism from $D$ to $M_s(\lambda)$ mapping~$v_D$ to~$v_M$.
\end{Lemma}

\begin{proof}
Since $t{\mathfrak h}[t]v_M = 0$ and the ${\mathfrak h}$-weights of both $v_D$ and $v_M$ are $\lambda$, it suf\/f\/ices to
show $x_\alpha^{\rho(\alpha)+1}v_M =0$ for all $\alpha \in \widehat{\Delta}^{\mathrm{re}}_+$.

First we consider the case $s=1$.
The assertion for $\alpha = -(\alpha_{p,n} + \alpha_{q,n-1}) + \delta$ with $p < q < n$ is proved by applying
Lemma~\ref{Lem:el} with $\beta= -\alpha_{p,m} + \delta$ and $\gamma = -\alpha_{q,m'}$.
The assertion for remaining $\alpha \in \widehat{\Delta}^{\mathrm{re}}_+$ is easily proved from the def\/ining relations
of $M_1(\lambda)$.

Next we consider the case $s=n$.
If $\sum\limits_{j=p}^{n-3} \lambda_j \ge \lambda_{n-1}$, the assertion for $\alpha = -\alpha_{p,n-1} + \delta$ is
proved by applying the lemma with $\beta = -\alpha_{p,n-2} + \delta$, $\gamma = -\alpha_{n-1}$.
Otherwise it is proved by applying the lemma with $\beta = -\alpha_{p,n-3}$, $\gamma = -\alpha_{n-2,n-1}+\delta$.
The assertion for $\alpha = -(\alpha_{p,n} + \alpha_{n-1}) + \delta$ is similarly proved, and then the assertion for
$\alpha = -(\alpha_{p,n} + \alpha_{q,n-1}) + \delta$ is verif\/ied by applying the lemma with $\beta = -(\alpha_{p,n} +
\alpha_{n-1}) + \delta$ and $\gamma = - \alpha_{q,n-2}$.
Finally the assertion for $\alpha = -(\alpha_{p,n} + \alpha_{q,n-1}) + 2\delta$ is shown by applying the lemma with
$\beta = -\alpha_{p,n} + \delta$ and $\gamma = -\alpha_{q,n-1} + \delta$.
The assertion for remaining $\alpha \in \widehat{\Delta}^{\mathrm{re}}_+$ is easily proved from the def\/ining relations
of~$M_n(\lambda)$.
\end{proof}

Now we can prove the existence of a~surjective ${\mathfrak g}[t]$-module homomorphism
\begin{gather*}
D\big(w_\circ \xi_n^{(s)}, w_\circ \xi_1^{(s)}, \ldots, w_\circ \xi_{n-1}^{(s)}\big) \twoheadrightarrow M_s(\lambda)
\end{gather*}
by exactly the same arguments with~\cite[two paragraphs below Lemma 5.2]{N} from Lemma~\ref{Lem:hom}.
Since $D\big(w_\circ \xi_n^{(s)}, w_\circ \xi_1^{(s)}, \ldots w_\circ \xi_{n-1}^{(s)}\big)\cong D\big(w_\circ
\xi_1^{(s)},\ldots,w_\circ \xi_n^{(s)}\big)$ holds by def\/inition, the proof is complete.

\subsection[Proof for $L(\pi_s) \twoheadrightarrow D\big(w_\circ \xi_1^{(s)},\ldots,w_\circ \xi_{n}^{(s)}\big)$]{Proof
for $\boldsymbol{L(\pi_s) \twoheadrightarrow D\big(w_\circ \xi_1^{(s)},\ldots,w_\circ \xi_{n}^{(s)}\big)}$}
\label{subsection:3-3}

We shall prove the assertion in the case $\varepsilon = +$ (the case $\varepsilon = -$ is similarly proved).
For $i \in I$, set
\begin{gather*}
p_i =
\begin{cases}
n - 2 -i, & i \le n-2,
\\
1, &  i=n-1,n,
\end{cases}
\end{gather*}
which is the distance between the nodes $i$ and $n-2$ in the Dynkin diagram.
We need the following lemma that is proved from~\cite{MR1883181}.

\begin{Lemma}\label{Lem:3-1}\quad
\begin{enumerate}\itemsep=0pt
\item[{\normalfont(i)}] Let $i_1, \ldots,i_p \in I$, $b_1,\ldots,b_p \in \mathbb{C}(q)^\times$ and $l_1, \ldots, l_p \in
\mathbb{Z}_{>0}$, and assume that
\begin{gather*}
b_r\notin q^{\mathbb{Z}_{> 0}} b_sq^{l_r-l_s+|p_{i_r} - p_{i_s}|}
\qquad
\text{for all}
\quad
r < s.
\end{gather*}
Then the submodule of $L_q\big(\boldsymbol{\pi}_{l_1,b_1}^{(i_1)}\big) \otimes \dots \otimes
L_q\big(\boldsymbol{\pi}_{l_p,b_p}^{(i_p)}\big)$ generated by the tensor product of highest $\ell$-weight vectors is isomorphic
to $L_q\Big(\prod\limits_{k =1}^p \boldsymbol{\pi}_{l_k,b_k}^{(i_k)}\Big)$.

\item[{\normalfont(ii)}] If $i,j \in I$, $b \in \mathbb{C}(q)^\times$, $l \in \mathbb{Z}_{>0}$ and $-|p_i - p_j| \le k \le |p_i
- p_j|$, then $L_q\big(\boldsymbol{\pi}_{l,b}^{(i)}\big) \otimes L_q\big(\boldsymbol{\pi}_{l,bq^k}^{(j)}\big)$ is simple.
\end{enumerate}
\end{Lemma}

\begin{proof}
(i) For $r < s$, it is directly checked that $L_q\big(\boldsymbol{\pi}_{l_s,b_s}^{(i_s)}\big) \otimes
L_q\big(\boldsymbol{\pi}_{l_r,b_r}^{(i_r)}\big)$ satisf\/ies the condition of~\cite[Corollary~6.2]{MR1883181}, which assures that
the module is generated by the tensor product of highest $\ell$-weight vectors.
Hence
\begin{gather*}
L_q\big(\boldsymbol{\pi}_{l_p,b_p}^{(i_p)}\big) \otimes L_q\big(\boldsymbol{\pi}_{l_{p-1},b_{p-1}}^{(i_{p-1})}\big) \otimes \cdots
\otimes L_q\big(\boldsymbol{\pi}_{l_1,b_1}^{(i_1)}\big)
\end{gather*}
is also generated by the tensor product of highest $\ell$-weight vectors (see~\cite[sentences above Corollary~6.2]{MR1883181}).
Now the assertion (i) follows by dualizing the statement and applying Lemma~\ref{Lem:dual} (the bijection $i \mapsto
\bar{i}$ in the lemma is $\tau_{n-1,n}$).
(ii) We see from the above argument that $L_q(\boldsymbol{\pi}_{l,b}^{(i)}) \otimes L_q(\boldsymbol{\pi}_{l,bq^k}^{(j)})$ is both cyclic
and cocyclic, and hence simple.
\end{proof}

\begin{Lemma}
\label{Lem:3-2}
Let $i \in I \setminus\{n-2,n-1\}$, $b \in \mathbb{C}^\times q^\mathbb{Z}$ and $l \in \mathbb{Z}_{\ge 0}$.
The graded limit $L\big(\boldsymbol{\pi}_{l,b}^{(i)}\boldsymbol{\pi}_{l,bq^{p_i-1}}^{(n-1)}\big)$ is isomorphic to
$D\big(lw_\circ(\varpi_i + \varpi_{n-1} + \Lambda_0)\big)$ as a~${\mathfrak g}[t]$-module.
\end{Lemma}

\begin{proof}
Note that
\begin{gather*}
\boldsymbol{\pi}_{l,b}^{(i)} \boldsymbol{\pi}_{l,bq^{p_i-1}}^{(n-1)} =
\begin{cases}
\boldsymbol{\pi}_{n}^{+}\big(l(\varpi_i + \varpi_{n-1}), bq^{(\delta_{i1}-1)l- i+1}\big) & \text{if} \quad i \le n-3,
\vspace{1mm}\\
\boldsymbol{\pi}_{1}^{+}\big(l(\varpi_{n} + \varpi_{n-1}), bq^{-l - n+2}\big) & \text{if} \quad i = n.
\end{cases}
\end{gather*}
Hence by Sections~\ref{subsection:3-1} and~\ref{subsection:3-2}, there exists a~surjective homomorphism
\begin{gather*}
D\big(lw_\circ(\varpi_i + \varpi_{n-1} + \Lambda_0)\big) \twoheadrightarrow
L\big(\boldsymbol{\pi}_{l,b}^{(i)}\boldsymbol{\pi}_{l,bq^{p_i-1}}^{(n-1)}\big).
\end{gather*}
Therefore it suf\/f\/ices to show the equality of the dimensions.
By Lemma~\ref{Lem:3-1}(ii), we have
\begin{gather*}
L_q\big(\boldsymbol{\pi}_{l,b}^{(i)}\boldsymbol{\pi}_{l,bq^{p_i-1}}^{(n-1)}\big) \cong L_q\big(\boldsymbol{\pi}_{l,b}^{(i)}\big) \otimes
L_q\big(\boldsymbol{\pi}_{l,bq^{p_i-1}}^{(n-1)}\big),
\end{gather*}
which implies
\begin{gather*}
\dim_\mathbb{C}L\big(\boldsymbol{\pi}_{l,b}^{(i)}\boldsymbol{\pi}_{l,bq^{p_i-1}}^{(n-1)}\big)
=\dim_{\mathbb{C}(q)}L_q\big(\boldsymbol{\pi}_{l,b}^{(i)}\boldsymbol{\pi}_{l,bq^{p_i-1}}^{(n-1)}\big)
\\
\qquad{}
=\dim_{\mathbb{C}(q)}L_q\big(\boldsymbol{\pi}_{l,b}^{(i)}\big)\cdot \dim_{\mathbb{C}(q)}
L_q\big(\boldsymbol{\pi}_{l,bq^{p_i-1}}^{(n-1)}\big) = \dim_{\mathbb{C}}L\big(\boldsymbol{\pi}_{l,b}^{(i)}\big)\cdot \dim_{\mathbb{C}}
L\big(\boldsymbol{\pi}_{l,bq^{p_i-1}}^{(n-1)}\big)
\\
\qquad{}
=\dim_{\mathbb{C}}D\big(lw_\circ(\varpi_i + \Lambda_0)\big) \cdot \dim_\mathbb{C} D\big(lw_\circ(\varpi_{n-1}+\Lambda_0)\big),
\end{gather*}
where the last equality follows from~\cite[Proposition~5.1.3]{MR2238884}.
On the other hand, we have
\begin{gather*}
\dim_\mathbb{C} D\big(lw_\circ(\varpi_i + \varpi_{n-1} + \Lambda_0)\big) = \dim_\mathbb{C}
D\big(lw_\circ(\varpi_i+\Lambda_0)\big) \cdot \dim_\mathbb{C} D\big(lw_\circ (\varpi_{n-1}+\Lambda_0)\big)
\end{gather*}
by~\cite[Theorem 1]{MR2235341}.
Hence the assertion is proved.
\end{proof}

Now let us begin the proof of the assertion $L(\boldsymbol{\pi}_s) \twoheadrightarrow
D\big(w_\circ\xi^{(s)}_1,\ldots,w_\circ\xi^{(s)}_{n}\big)$.
First we prove this in the case $s=n$.
Let $(a_i)_{i \in I}$ be the sequence in the def\/inition of $\boldsymbol{\pi}_n= \boldsymbol{\pi}_n^{+}(\lambda,a)$, and def\/ine
$U_q(\mathbf{L}{\mathfrak g})$-modules $L_{q}[j]$ for $1 \le j \le n-2$ and $j = n$ by
\begin{gather}
\label{eq:Lqj}
L_q[j] =
\begin{cases}
L_q\big(\boldsymbol{\pi}_{\lambda_j,a_j}^{(j)}\big), & \text{if} \quad j<\ell, \quad j=n-2,n,
\vspace{1mm}\\
L_q\big(\boldsymbol{\pi}_{\lambda_\ell - \bar{\lambda},a_\ell q^{-\bar{\lambda}}}^{(\ell)}\big)\otimes
L_q\big(\boldsymbol{\pi}_{\bar{\lambda},a_\ell q^{\lambda_\ell-\bar{\lambda}}}^{(\ell)} \boldsymbol{\pi}_{\bar{\lambda}, a_\ell
q^{\lambda_\ell - \bar{\lambda} + p_\ell- 1}}^{(n-1)}\big), & \text{if} \quad j = \ell,
\vspace{1mm}\\
L_q\big(\boldsymbol{\pi}_{\bar{\lambda}, a_1q^{-\lambda_1-\bar{\lambda}+n-4}}^{(n-1)}\big) \otimes
L_q\big(\boldsymbol{\pi}_{\lambda_1,a_1}^{(1)}\boldsymbol{\pi}_{\lambda_1, a_1q^{ n - 4}}^{(n-1)} \big), & \text{if} \quad \ell=0, \quad j=1,
\vspace{1mm}\\
L_q\big(\boldsymbol{\pi}_{\lambda_j,a_j}^{(j)}\boldsymbol{\pi}_{\lambda_j,a_jq^{p_j-1}}^{(n-1)}\big), & \text{otherwise}.
\end{cases}\!\!\!\!\!\!\!
\end{gather}
There exists an injective $U_q(\mathbf{L}{\mathfrak g})$-module homomorphism
\begin{gather*}
L_q(\boldsymbol{\pi}_n) \hookrightarrow L_q[1] \otimes \dots \otimes L_q[n-2] \otimes L_q[n]
\end{gather*}
by Lemma~\ref{Lem:3-1}, and this induces a~$U_{\bf{A}}(\mathbf{L}{\mathfrak g})$-module homomorphism
\begin{gather*}
L_{\bf{A}}(\boldsymbol{\pi}_n) \to L_{\bf{A}}[1] \otimes \dots \otimes L_{\bf{A}}[n-2] \otimes L_{\bf{A}}[n],
\end{gather*}
where we set
\begin{gather*}
L_{\bf{A}}[j] =
\begin{cases}
L_{\bf{A}}(\boldsymbol{\pi}), & \text{if} \quad L_q[j]= L_q(\boldsymbol{\pi}),
\\
L_{\bf{A}}\big(\boldsymbol{\pi}^1\big) \otimes L_{\bf{A}}\big(\boldsymbol{\pi}^2\big), & \text{if} \quad L_q[j] = L_q\big(\boldsymbol{\pi}^1\big) \otimes L_q\big(\boldsymbol{\pi}^2\big).
\end{cases}
\end{gather*}
Applying $\mathbb{C} \otimes_{\bf A} -$ and taking the pull-back with respect to the automorphism $\varphi_{\bar{a}}$,
we obtain a~${\mathfrak g}[t]$-module homomorphism $L(\boldsymbol{\pi}_n) \to \bigotimes_{j} L[j]$, where $L[j]$ denotes the
graded limit or the tensor product of the two graded limits.
Note that, by construction, this homomorphism maps a~highest weight vector of $L(\boldsymbol{\pi}_n)$ to the tensor product of
highest weight vectors.
By Lemma~\ref{Lem:3-2} and~\cite[Proposition~5.1.3]{MR2238884}, we have
\begin{gather*}
L[j] \cong
\begin{cases}
D\big((\lambda_\ell-\bar{\lambda})w_\circ(\varpi_\ell + \Lambda_0)) \otimes D\big(\bar{\lambda} w_\circ(\varpi_\ell +
\varpi_{n-1} + \Lambda_0)\big), & \text{if} \quad j = \ell,
\\
D\big(\bar{\lambda}w_\circ(\varpi_{n-1}+\Lambda_0)\big) \otimes D\big(\lambda_1w_\circ(\varpi_1 + \varpi_{n-1} +
\Lambda_0)\big), & \text{if} \quad \ell =0, \quad j =1,
\\
D\big(w_\circ \xi_{j}^{(n)}\big), & \text{otherwise}.
\end{cases}
\end{gather*}
Hence in order to complete the proof, it suf\/f\/ices to show that
\begin{gather*}
\begin{split}
& D\big((\lambda_\ell-\bar{\lambda})w_\circ(\varpi_\ell + \Lambda_0),\bar{\lambda} w_\circ(\varpi_\ell + \varpi_{n-1} +
\Lambda_0)\big) \cong D\big(w_\circ \xi_\ell^{(n)}\big),
\qquad
\text{and}
\\
& D\big(\bar{\lambda}w_\circ(\varpi_{n-1}+\Lambda_0), \lambda_1w_\circ(\varpi_1 + \varpi_{n-1} + \Lambda_0)\big) \cong
D\big(w_\circ \xi_1^{(n)}\big)
\qquad
\text{when}
\quad
\ell = 0.
\end{split}
\end{gather*}
The f\/irst isomorphism follows since we have
\begin{gather*}
D\big(w_\circ \xi_\ell^{(n)}\big) \cong F_{w_\circ \sigma^\ell} \big(\Lambda_\ell^{(n)}\big) \cong F_{w_\circ
\sigma^\ell} \big(D\big((\lambda_\ell - \bar{\lambda})\Lambda_0\big) \otimes D\big(\bar{\lambda}(\varpi_{n-1} +
\Lambda_0)\big)\big)
\\
\hphantom{D\big(w_\circ \xi_\ell^{(n)}\big)}{}
\cong D\big((\lambda_\ell-\bar{\lambda})w_\circ(\varpi_\ell + \Lambda_0), \bar{\lambda} w_\circ(\varpi_\ell +
\varpi_{n-1} + \Lambda_0)\big)
\end{gather*}
by Proposition~\ref{Prop:construction}, and the second one is also proved similarly.
The proof for $s=n$ is complete.

The case $s=1$ can be proved by a~similar (and simpler) argument in which we replace the def\/inition of~$L_q[i]$ given
in~\eqref{eq:Lqj} by the following:
\begin{gather*}
L_q[j] =
\begin{cases}
L_q\big(\boldsymbol{\pi}_{\lambda_j,a_j}^{(j)}\big), & j \le n-2,
\vspace{1mm}\\
L_q\big(\boldsymbol{\pi}_{\lambda_{m'},a_{m'}}^{(m')}\boldsymbol{\pi}_{\lambda_{m'},a_{m'}}^{(m)}\big), & j = n-1,
\vspace{1mm}\\
L_q\big(\boldsymbol{\pi}_{\lambda_{m}-\lambda_{m'},a_{m'}q^{\lambda_m}}^{(m)}\big), & j = n.
\end{cases}
\end{gather*}

\subsection*{Acknowledgements}

The author would like to thank Steven V.~Sam for informing him of the results in~\cite{Sam}.
This work was supported by JSPS Grant-in-Aid for Young Scientists (B) No.~25800006, and by World Premier International
Research Center Initiative (WPI Initiative), MEXT, Japan.

\pdfbookmark[1]{References}{ref}
\LastPageEnding

\end{document}